\documentclass[a4paper,10pt]{amsart}
\usepackage{amsmath}
\usepackage{amsthm}

\usepackage[all]{xy}

\theoremstyle{plain}
\newtheorem{thm}{Theorem}
\newtheorem{prop}[thm]{Proposition}
\newtheorem{lem}[thm]{Lemma}
\newtheorem{cor}[thm]{Corollary}
\theoremstyle{remark}

\newcommand{\acat}[1]{\mathbf{a}_{#1}}
\newcommand{\Acat}[1]{\mathbf{A}_{#1}}
\newcommand{\aff}{\mathbf{A}}
\newcommand{\ccat}{\mathbf{C}}
\newcommand{\cacat}[1]{\mathbf{\hat{a}}_{#1}}
\newcommand{\cc}{\mathbf{C}}
\newcommand{\cofilt}{\mathbf{CoFilt}(\mfam)}
\newcommand{\defm}{\textup{Def}_{\mfam}}
\newcommand{\extn}{\mathbf{Ext}(\mfam)}
\newcommand{\extncl}{\mathbf{E}}
\newcommand{\filt}{\mathbf{Filt}(\mfam)}
\newcommand{\Gm}{\Gamma}
\newcommand{\grdmod}{\mathbf{gr}\modc(D)}
\newcommand{\grhol}{\mathbf{grHol}}
\newcommand{\indextncl}{\mathbf{IE}}
\newcommand{\indmodfcl}{\mathbf{IM}}
\newcommand{\infdef}{\mathbf{T}^1}
\newcommand{\mfam}{\mathcal{F}}
\newcommand{\modc}{\mathbf{Mod}}
\newcommand{\modf}{\mathbf{Mod}(\mfam)}
\newcommand{\modfcl}{\mathbf{M}}
\newcommand{\modw}{\mathbf{fdMod}(W)}
\newcommand{\n}{\mathbf{N}}
\newcommand{\nsextncl}{\mathbf{*E}}
\newcommand{\nsmodfcl}{\mathbf{*M}}
\newcommand{\obstrdef}{\mathbf{T}^2}
\newcommand{\op}[1]{{#1}^{\text{op}}}
\newcommand{\rmod}{\modc(R)}
\newcommand{\sets}{\mathbf{Sets}}
\newcommand{\subfam}{\mfam'}
\newcommand{\vect}[1]{\underline{#1}}
\newcommand{\weyl}{A_1(k)}
\newcommand{\weylf}{B_1(k)}
\newcommand{\x}{\mathbf{X}(\mfam,\Gm)}
\newcommand{\z}{\mathbf{Z}}

\DeclareMathOperator{\aut}{Aut}

\DeclareMathOperator{\der}{Der}
\DeclareMathOperator{\diff}{Diff}
\DeclareMathOperator{\enm}{End}
\DeclareMathOperator{\ext}{Ext}
\DeclareMathOperator{\gr}{gr}
\DeclareMathOperator{\hmm}{Hom}
\DeclareMathOperator{\id}{id}
\DeclareMathOperator{\mor}{Mor}
\DeclareMathOperator*{\osum}{\oplus}

\DeclareMathOperator{\re}{Re}
\DeclareMathOperator{\soc}{soc}

\begin{document}

\title{Iterated extensions in module categories}
\author{Eivind Eriksen}
\address{University of Oslo, P.O. Box 1053 Blindern, 0316 Oslo,
Norway}
\email{eeriksen@math.uio.no}

\thanks{This research has been supported by a Marie Curie Fellowship 
of the European Community programme "Improving Human Research Potential
and the Socio-economic Knowledge Base" under contract number 
HPMF-CT-2000-01099}

\begin{abstract}
Let $k$ be an algebraically closed field, let $R$ be an associative $k$-algebra,
and let $\mfam = \{ M_{\alpha}: \alpha \in I \}$ be a family of orthogonal 
points in $\rmod$ such that $\enm_R(M_{\alpha}) \cong k$ for all $\alpha \in I$. 
Then $\modf$, the minimal full sub-category of $\rmod$ which contains $\mfam$ 
and is closed under extensions, is a full exact Abelian sub-category of $\rmod$ 
and a length category in the sense of Gabriel \cite{ga73}. 

In this paper, we use iterated extensions to relate the length category $\modf$
to noncommutative deformations of modules, and use some new methods to study  
$\modf$ via iterated extensions. In particular, we give a new proof of the 
characterization of uniserial length categories, which is constructive. As an 
application, we give an explicit description of some categories of holonomic and 
regular holonomic $D$-modules on curves which are uniserial length categories.  
\end{abstract}

\maketitle

\section*{Introduction}

Let $\ccat = \rmod$ be the category of left modules over an 
associate ring $R$. We shall assume that $\mfam = \{ M_{\alpha}:
\alpha \in I \}$ is a family of non-zero, pairwise non-isomorphic
objects in $\ccat$, and consider the minimal full sub-category 
$\modf$ of $\ccat$ which contains $\mfam$ and is closed under 
extensions. 

An alternative and explicit description of $\modf$ is useful: An 
object $M$ of $\ccat$ is in $\modf$ if and only if there is a 
cofiltration
\[ M = C_n \to C_{n-1} \to \dots \to C_1 \to C_0 = 0 \]
in $\ccat$ such that $f_i: C_i \to C_{i-1}$ is surjective with 
kernel $K_i \cong  M_{\alpha(i)}$ with $\alpha(i) \in I$ for $1 
\le i \le n$. Equivalently, $M$ is in $\modf$ if and only if 
there is a filtration 
\[ 0 = F_n \subseteq F_{n-1} \subseteq \dots \subseteq F_1 
\subseteq F_0 = M \]
in $\ccat$ such that $K_i = F_{i-1}/F_i \cong M_{\alpha(i)}$ with 
$\alpha(i) \in I$ for $1 \le i \le n$. 

In general, $\modf$ is not an exact Abelian sub-category of $\ccat$: 
It does not necessarily contain its kernels and cokernels. Ringel 
\cite{ri76} has shown that $\modf$ is a full, extension closed and 
exact Abelian sub-category of $\ccat$ with $\mfam$ as its simple 
objects if and only if $\mfam$ is a family of orthogonal points in 
$\ccat$. In this paper, we shall assume that this is the case. So by
definition, $\enm_R(M_{\alpha})$ is a division ring and 
$\hmm_R(M_{\alpha},M_{\beta})=0$ for all $\alpha, \beta \in I$. This 
means that $\modf$ is a length category in the sense of Gabriel 
\cite{ga73}. 

Let us assume that the family $\mfam$ of orthogonal points in $\ccat$ 
is given. We shall consider the following problem: Classify the 
indecomposable objects in $\modf$, up to isomorphism. This problem is
fundamental, but there is nothing original about it -  we cite Gabriel 
\cite{ga73}:
\begin{quote}
The main and perhaps hopeless purpose of representation theory is to 
find an efficient general method for constructing the indecomposable 
objects by means of the simple objects, which are supposed to be 
given.
\end{quote} 
It is not plausible to expect a full solution to this problem. It is 
well-known that the problem is \emph{wild} in many cases, for instance
when $\mfam$ is a complete family of simple left modules over the first
Weyl algebra $\weyl$ over any algebraically closed field $k$ of 
characteristic $0$.

It is maybe better to consider the following problem: Give necessary
and sufficient conditions for the category $\modf$ to be \emph{tame}.
We would of course also like to classify the indecomposable modules in
$\modf$ in these cases. We are not able to give a complete solution to 
this problem at present. However, we shall introduce some new methods 
which we believe are useful for attacking the problem, and we shall give 
the solution to the problem in some special cases. 

Before we go on to study the above problem in more detail, we remark 
that so far, we have only treated the case when $\ccat = \rmod$ is the 
category of left modules over an associative ring $R$. It makes sense 
to consider the problem above when $\ccat$ is any Abelian category. 
Most of the results will remain true with mild restrictions (and in 
many cases, none) on the Abelian category $\ccat$. We shall only treat 
the case $\ccat = \rmod$ in this paper, since this will lead to a much 
more readable exposition, and leave it to the reader to figure out how 
to generalize the results to more general Abelian categories $\ccat$. 
The only exception is one of our applications, where we assume that 
$\ccat$ is the category of graded modules over a graded ring $R$, in 
which case all results of this paper remain valid. 

Let us consider the category $\extn$ of \emph{iterated extensions}
of the family $\mfam$: The objects of $\extn$ are couples $(M,C)$, where
$M$ is an object of $\ccat$ and $C$ is a cofiltration of $M$ of the type
considered above, see section \ref{s:extn} for details. Clearly, there
is a forgetful functor $\extn \to \rmod$, and image of this functor is 
exactly the category $\modf$. The reason why we would like to consider
the category $\extn$, is that many useful invariants are naturally 
defined there.

Let $(M,C)$ be an object of $\extn$. We define the length $n$ to be the 
length $n$ of the cofiltration $C$, and the order vector $\vect \alpha
\in I^n$ to be the vector defined by $K_i \cong M_{\alpha(i)}$ for $1 
\le i \le n$. All objects in an isomorphism class of $\extn$ has the 
same length and order vector. We define the \emph{extension type} of 
the object $(M,C)$ to be the ordered quiver $\Gm$ with vertices 
$\{ \alpha(i) \in I: 1 \le i \le n \}$, and arrows $\{ a_{i-1,i}: 2 \le 
i \le n \}$, where $a_{12} < a_{23} < \dots < a_{n-1,n}$ is the total 
ordering of the arrows, and $a_{i-1,i}$ is an arrow from vertice 
$\alpha(i-1)$ to vertice $\alpha(i)$. The extension type $\Gm$ is a 
convenient way of representing the invariants given by the length $n$ and 
the order vector $\vect \alpha$. 

Let us assume that $R$ is an algebra over an algebraically closed 
(commutative) field $k$. Under some finiteness conditions, we shall show 
that the category $\extn$, and therefore also the category $\modf$, is 
determined by the noncommutative deformations of the family $\mfam$.  

Let $(M,C)$ be an object of $\extn$, and let $\Gm$ be the extension type 
of $(M,C)$. We associate with $\Gm$ the $k$-algebra $k[\Gm]$, see section
\ref{s:extype} for details. This is a $p$-pointed $k$-algebra, where $p$ 
is the number of vertices in $\Gm$, so there are natural maps $k^p \to 
k[\Gm] \to k^p$. Moreover, the radical $I$ of $k[\Gm]$ satisfies $I^n = 0$, 
where $n$ is the length of $(M,C)$. So $k[\Gm]$ is a complete Artinian ring 
in the $I$-adic topology, and therefore $k[\Gm]$ is an object in the 
category $\acat p(n)$ of complete Artinian $p$-pointed algebras.  

There is a theory of noncommutative deformations of modules, due to 
Laudal, see Laudal \cite{lau95}, \cite{lau98}, \cite{lau00}. We refer to 
the preprint Eriksen \cite{er02} for a convenient form of the results we 
need in this paper. Let us recall the main theorem: For any finite family 
$\mfam = \{ M_{\alpha}: \alpha \in I \}$ of left $R$-modules, there is a 
noncommutative deformation functor $\defm: \acat p \to \sets$. If 
$\ext^i(M_{\alpha},M_{\beta})$ is a finite dimensional vector space over 
$k$ for $i=1,2$ and for all $\alpha, \beta \in I$, then $\defm$ has a 
pro-representing hull $H$, which is unique up to (non-canonical) isomorphism 
in $\cacat p$. Both the hull $H$ and the corresponding versal family $M_H 
\in \defm(H)$ are in principle constructible.

Let $\Gm$ be a fixed extension type, and consider iterated extensions
of the family $\mfam$ with extension type $\Gm$. Clearly, these are all 
iterated extensions of the finite sub-family of $\mfam$ given by the
vertices of $\Gm$. Let us denote by $\extncl(\mfam,\Gm)$ the set of 
isomorphism classes of iterated extensions in $\extn$ with extension 
type $\Gm$. We give a proof of the following result, due to Laudal:

\begin{thm}[Laudal]
Let $\mfam$ be a finite family of $p$ orthogonal points in $\rmod$, and let 
$\Gm$ be the extension type of some object $(M,C)$ of $\extn$ with $p$ 
vertices. Then there is a bijective correspondence between $\defm(k[\Gm])$ 
and $\extncl(\mfam,\Gm)$.    
\end{thm}

We show that $\x = \mor(H,k[\Gm])$ has a natural structure of an affine 
scheme over $k$, and by definition of the pro-representing hull $H$, there 
is a surjection $\x \to \extncl(\mfam,\Gm)$. The forgetful functor $\extn 
\to \modf$ maps the isomorphism classes of $\extncl(\mfam,\Gm)$ to a subset 
of isomorphism classes of modules in $\modf$, and we shall denote this set 
of isomorphism classes by $\modfcl(\mfam,\Gm)$. So there are natural 
surjections $\x \to \extncl(\mfam,\Gm) \to \modfcl(\mfam,\Gm)$. In 
particular, the sets $\extncl(\mfam,\Gm)$ and $\modfcl(\mfam,\Gm)$ are 
quotients of the affine variety $\x$. 

The \emph{species} of $\modf$ is given by $(K_{\alpha},E_{\alpha,\beta})$, 
where $K_{\alpha} = \enm_R(M_{\alpha})$ is a division ring and 
$E_{\alpha,\beta} = \ext^1_R(M_{\alpha},M_{\beta})$ is a 
$K_{\beta}$-$K_{\alpha}$ bimodule for all $\alpha, \beta \in I$. From now 
on, we shall assume that $R$ is an algebra over an algebraically closed 
field $k$, that $\mfam$ is a family of orthogonal points in $\rmod$, and 
that $\enm_R(M_{\alpha}) \cong k$ for all $\alpha \in I$. In this case, the
species of $\modf$ is called a $k$-quiver, because it is completely 
determined by the \emph{Gabriel quiver}, defined by the set of vertices $I$ 
and $\dim_k E_{\alpha,\beta}$ arrows from $\alpha$ to $\beta$ for each pair 
of vertices $\alpha, \beta \in I$. 

It is known that the species, and therefore the Gabriel quiver, contains a 
lot of information about the length category $\modf$. In fact, we see from 
Laudals theorem above that the only information which is not present in the 
Gabriel quiver is the obstruction theory of the family $\mfam$. Under the
same finiteness conditions as in Laudals theorem, it is shown in Deng, Xiao 
\cite{de-xi98} that if $\modf$ is a hereditary category, then it is 
equivalent to the category of \emph{small representations} of the Gabriel 
quiver of $\mfam$. A representation is called small if it is nilpotent (and 
finite dimensional). 

The hereditary case mentioned above is an example of an \emph{unobstructed 
case}, where $\modf$ is completely determined by its species. If the Gabriel
quiver is without loops, the category of small representations of the quiver
is just the usual category of finite dimensional representations of the 
quiver. So in this case, it is well-known when the category $\modf$ is wild,
tame and finite. If the Gabriel quiver has loops, we ask when the category of 
its small representations is wild, tame and finite. We do not know the 
complete answer to this question. 

We say that a module $M$ in $\modf$ is \emph{uniserial} if the lattice of 
submodules of $M$ is a chain. Moreover, we say that the length category 
$\modf$ is \emph{uniserial} if all indecomposable modules in $\modf$ are 
uniserial. It is maybe not so obvious, but the uniserial case is also 
unobstructed in the sense that the obstructions for deforming the family 
$\mfam$ do not survive in the category $\modf$. In other words, the 
category $\modf$ is completely determined by its species in the uniserial 
case.

\begin{thm}
Let $\mfam$ be a family of orthogonal points such that $\enm_R(M_{\alpha}) 
\cong k$ for all $\alpha \in I$. Then the category $\modf$ is uniserial if
and only if each connected component of the Gabriel quiver of $\mfam$ is 
either a cycle or a linear quiver. Moreover, if this is the case, then 
$\indmodfcl(\mfam,\Gm)$ has a single element if $\Gm$ is admissible, and 
otherwise it is empty.
\end{thm}

By definition, an extension type $\Gm$ is \emph{admissible} if the 
corresponding path appears in the Gabriel quiver of $\modf$. The set  
$\indmodfcl(\mfam,\Gm)$ denotes the subset of $\modfcl(\mfam,\Gm)$ 
consisting of indecomposable isomorphism classes.
 
The first part of this theorem is known, see Gabriel \cite{ga73}, Amdal and
Ringdal \cite{am-ri68}. As far as we know, our proof is new. It has the good
property of being constructive, which is manifested in the second part of the 
theorem. In fact, we can construct the indecomposable modules corresponding
to admissible extension types explicitly. 

We take advantage of this fact when we apply the theorem to some categories 
of regular holonomic $D$-modules over curves of characteristic $0$. In 
particular, we show that the category of graded holonomic $D$-modules over 
the first Weyl algebra is uniserial when $k$ has characteristic $0$. We 
also describe the graded holonomic modules explicitly in this case.

Clearly, the length category $\modf$ is tame (or even finite) when the 
condition of the theorem holds (that is, when $\modf$ is uniserial). On the 
other hand, a sufficient condition for the category $\modf$ to be wild in a 
strong sense is known. To be more precise, let $W = k <x,y>$ be the free 
associative $k$-algebra on two generators, and let $\modw$ be the category of 
left $W$-modules which are finite dimensional as vector spaces over $k$. We 
shall say that the length category $\modf$ is \emph{wild} if there is a full 
exact embedding of $\modw$ into $\modf$. 

We recall a well-known argument to show how hopeless it is to classify the
indecomposable modules when $\modf$ is wild: Let $R$ be any $k$-algebra
which has finite dimension as vector space over $k$. Then there is a full
exact embedding of the category of left $R$-modules which have finite dimension
over as vector spaces over $k$ into $\modw$, and therefore into $\modf$.
Since full exact embeddings preserve indecomposable modules, a classification
of indecomposable modules in $\modf$ would contain a classification of all
finite dimensional indecomposable modules over all finite dimensional 
$k$-algebras.

The following theorem, essentially due to Klingler, Levy \cite{kl-le95}, gives 
a sufficient condtion for the length category $\modf$ to be wild in the above 
sense:

\begin{thm}[Klingler-Levy]
Let $\mfam$ be a family of orthogonal points such that $\enm_R(M_{\alpha}) 
\cong k$ for all $\alpha \in I$. If the Gabriel quiver of $\modf$ contains 
the quiver $Q_5$ given by 
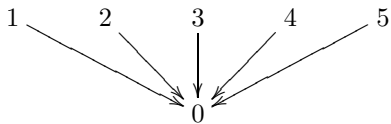
\[ \xymatrix{ 1 \ar[drr] & 2 \ar[dr] & 3 \ar[d] & 4 \ar[dl] & 5 \ar[dll] \\
& & 0 & & } \] 
then there is a full exact embedding of $\modw$ into $\modf$, and all modules 
in the image of this embedding has socle-height $2$. In particular, $\modf$ is 
wild in this case.
\end{thm}

We have described two extreme cases: The uniserial case, where we have found a
very nice description of the length category $\modf$, and a wild case, where it
would be hopeless to describe $\modf$. It would be interesting to know what 
happens with $\modf$ in all the intermediate cases. We would expect that an
understanding of the obstructions of the deformations of $\mfam$ is necessary 
to give the full answer to this question. But as we have noted above, even in
the hereditary case this question is open.

\section{Categories of iterated extensions}
\label{s:extn}

Let $R$ be an associative ring, and let $\mfam = \{ M_{\alpha} :
\alpha \in I \}$ be a fixed family of non-zero, pairwise
non-isomorphic left $R$-modules. We shall define the category
$\extn$ of \emph{iterated extensions} of the family $\mfam$.

Let us first consider the category $\cofilt$ of \emph{modules with
cofiltration} over the family $\mfam$, defined in the following way:
An object of $\cofilt$ is a couple $(M,C)$, where $M$ is a left
$R$-module and $C$ is a cofiltration of $M$ of the form
    \[ M = C_n \to C_{n-1} \to \dots \to C_2 \to C_1 \to C_0 = 0, \]
such that $f_i: C_i \to C_{i-1}$ is surjective, and $K_i = \ker(f_i)
\cong M_{\alpha(i)}$ with $\alpha(i) \in I$ for $1 \le i \le n$. The
integer $n \ge 0$ is called the \emph{length} of the cofiltration,
and the modules $K_1, \dots, K_n$ are called the \emph{factors} of
the cofiltration. Let $(M,C)$ and $(M',C')$ be objects of $\cofilt$ of
lengths $n,n' \ge 0$, and let $N = \max \{ n,n' \}$. A morphism $\phi:
(M,C) \to (M',C')$ is a collection $\{ \phi_i \in \hmm_R(C_i,C'_i): 0
\le i \le N \}$ of $R$-linear homomorphisms such that $\phi_{i-1} f_i =
f'_i \phi_i$ for $1 \le i \le N$. By convention, $C_i=M$ if $i>n$ and
$C'_i=M'$ if $i>n'$.

Similarly, we consider the category $\filt$ of \emph{modules with
filtration} over the family $\mfam$, defined in the following way: An
object of $\filt$ is a couple $(M,F)$, where $M$ is a left $R$-module
and $F$ is a filtration of $M$ of the form
    \[ 0 = F_n \subseteq F_{n-1} \subseteq \dots \subseteq F_0 = M, \]
such that $K_i = F_{i-1}/F_i \cong M_{\alpha(i)}$ with $\alpha(i) \in
I$ for $1 \le i \le n$. The integer $n \ge 0$ is called the
\emph{length} of the filtration, and the modules $K_1, \dots, K_n$ are
called the \emph{factors} of the filtration. Let $(M,F)$ and $(M',F')$
be objects in $\filt$ with lengths $n,n' \ge 0$, and let $N = \max \{
n,n' \}$. A morphism $\phi: (M,F) \to (M',F')$ is a homomorphism $\phi
\in \hmm_R(M,M')$ such that $\phi(F_i) \subseteq F'_i$ for $1 \le i \le
N$. By convention, $F_i=0$ if $i>n$ and $F'_i=0$ if $i>n'$.

Clearly, the categories $\cofilt$ and $\filt$ are equivalent, since
filtrations and cofiltrations are dual notions: If a cofiltration $C$
of $M$ of length $n \ge 0$ is given, let $F$ be the filtration defined
by $F_i = \ker(M \to C_i)$ for $0 \le i \le n$. Then the assignment
$(M,C) \mapsto (M,F)$ defines a functor $\cofilt \to \filt$. Conversely,
if a filtration $F$ of $M$ of length $n \ge 0$ is given, let $C$ be the
cofiltration defined by $C_i = M/F_i$ for $0 \le i \le n$, with the
natural surjections $f_i: C_i \to C_{i-1}$. Then the assignment $(M,F)
\mapsto (M,C)$ defines a functor $\filt \to \cofilt$. We see that these
functors are inverses of each other, and therefore define an equivalence
of categories between $\cofilt$ and $\filt$. Moreover, this equivalence
preserves the length $n$ and the factors $K_1, \dots, K_n$.

We say that an object $(M,C)$ in $\cofilt$ is an \emph{iterated 
extension} of the family $\mfam$, and we define the category $\extn$ of
\emph{iterated extensions} of the family $\mfam$ to equal the
category $\cofilt$. Moreover, we say that the length of an iterated 
extension $(M,C)$ is the length $n$ of the cofiltration $C$, and the
factors of $(M,C)$ are the factors $K_1, \dots, K_n$ of the cofiltration
$C$.

Clearly, an iterated extension $(M,C)$ of the family $\mfam$ of
length $n \le 1$ is given in the following way: If $n=0$ then $M=0$,
with the trivial filtration $C_0=0$. If $n=1$, then $M \in \mfam$,
and the filtration $C$ is given by $C_1=M, \; C_0=0$.

As the name suggests, iterated extensions of the family $\mfam$
of length $n \ge 2$ can be characterized in terms of extensions. Recall
that given a pair $M',M''$ of left $R$-modules, $M$ is said to be an
\emph{extension} of $M''$ by $M'$ if there exists an exact sequence
$0 \to M' \to M \to M'' \to 0$ of left $R$-modules.

\begin{lem}
Let $M$ be a left $R$-module and let $n \ge 2$ be an integer. Then the
following conditions are equivalent:
\begin{enumerate}
\item There exists a cofiltration $C$ of $M$ such that $(M,C)$ is an
object of $\extn$ of length $n$,
\item $M$ is an extension of $M''$ by $M'$, where $(M',C'), (M'',C'')$
are objects of $\extn$ of lengths $n',n''$, with $n'+n''=n$ and $n',
n'' < n$.
\end{enumerate}
\end{lem}
\begin{proof}
If $(M,C)$ is an iterated extension of the family $\mfam$ of
length $n$, then $M$ is an extension of $C_{n-1}$ by $K_n$. But $K_n
\in \mfam$ and $C_{n-1}$ is clearly an iterated extension of the
family $\mfam$ of length $n-1$. For the other implication, assume that
$(M',C')$ and $(M'',C'')$ are iterated extensions of the family
$\mfam$ of lengths $n',n''$. We construct a cofiltration $C$ of $M$ of
length $n=n'+n''$ in the following way: Let $f:M' \to M$ and $g:M \to
M''$ be the maps given by the extension $0 \to M' \to M \to M'' \to 0$,
let $F'$ be the filtration of $M'$ corresponding to the cofiltration
$C'$, and let $F''$ be the filtration of $M''$ corresponding to $C''$.
We define $F_i = g^{-1}(F''_i)$ for $0 \le i \le n''$, and $F_i =
f(F'_{i-n''})$ for $n'' \le i \le n$. Then $F$ is a filtration of $M$,
and we have $F_{i-1}/F_i \cong \ker(M \to C''_{i-1})/\ker(M \to C''_i)
\cong K''_i$ for $0 \le i \le n''$, $F_{i-1}/F_i \cong K'_{i-n''}$
for $n'' \le i \le n$. Let $C$ be the cofiltration of $M$ corresponding
to the filtration $F$. Then $(M,C)$ is an iterated extension of
the family $\mfam$ of length $n=n'+n''$.
\end{proof}

Let $(M,C)$ be an iterated extension of the family $\mfam$ of
length $n$. For $1 \le i \le n$, we have $K_i = \ker(C_i \to C_{i-1})
\cong M_{\alpha(i)}$ for a unique $\alpha(i) \in I$. The resulting 
vector $\vect \alpha = (\alpha(1), \dots, \alpha(n)) \in I^n$ is called 
the \emph{order vector} of $(M,C)$. Clearly, it is uniquely defined by 
the cofiltration $C$ since the family $\mfam$ consists of pairwise
non-isomorphic modules.

Moreover, for $2 \le i \le n$, the filtration $C$ induces the
following commutative diagram of left $R$-modules
        \[
        \xymatrix{
        0 \ar[r] & K_i \ar[r] & C_i \ar[r]^{f_i} & C_{i-1} \ar[r] & 0 \\
        0 \ar[r] & K_i \ar[r] \ar@{=}[u] & \ker(f_{i-1} \circ f_i)
        \ar[r]_{f_i} \ar[u] & K_{i-1} \ar[r] \ar[u] & 0, \\
        }
        \]
where the rows are exact. We denote the extensions corresponding to the
upper and lower row by $\xi_i \in \ext^1_R(C_{i-1},K_i)$ and $\tau_i \in
\ext^1_R(K_{i-1},K_i)$ respectively. The commutativity of the above
diagram means that $\xi_i \mapsto \tau_i$ under the map
    \[ \ext^1_R(K_{i-1} \to C_{i-1},K_i): \ext^1_R(C_{i-1},K_i) \to
    \ext^1_R(K_{i-1},K_i) \]
induced by the inclusion $K_{i-1} \to C_{i-1}$.

Let $\phi: (M,C) \to (M',C')$ be an isomorphism in $\extn$. Then the
homomorphisms $\phi_i: C_i \to C'_i$ are all isomorphisms. This proves
the following useful result:

\begin{lem} \label{l:extnchar}
Let $\phi: (M,C) \to (M',C')$ be an isomorphism in $\extn$. Then
$(M,C)$ and $(M',C')$ have the same length and order vector. Moreover,
the isomorphism $\phi$ induces isomorphisms of Abelian groups
\begin{align*}
\ext^1_R(C_{i-1},K_i) &\to \ext^1_R(C'_{i-1},K'_i) \\
\ext^1_R(K_{i-1},K_i) &\to \ext^1_R(K'_{i-1},K'_i)
\end{align*}
for $2 \le i \le n$, where $n$ is the common length of $(M,C)$ and
$(M',C')$. Under these isomorphisms, $\xi_i \mapsto \xi'_i$ and $\tau_i
\mapsto \tau'_i$.
\end{lem}

There is a forgetful functor $\extn \to \rmod$ given by $(M,C) \mapsto M$,
where $\rmod$ denotes the category of left $R$-modules. The image of this
functor defines a full subcategory of $\rmod$, which we denote by $\modf$.
Clearly, a left $R$-module $M$ is an object of $\modf$ if and only if there
exists a cofiltration $C$ of $M$ such that $(M,C)$ is an object of $\extn$.
Lemma \ref{l:extnchar} translates to the following characterization of the
category $\modf$:

\begin{prop} \label{p:modfchar}
Let $\mfam$ be a family of non-zero, pairwise non-isomorphic left
$R$-modules. Then the category $\modf$ is the minimal full sub-category of
$\rmod$ which contains $\mfam$ and is closed under extensions.
\end{prop}

It follows that $\modf \subseteq \rmod$ is a full, exact subcategory which
is closed under extensions. But in general, $\modf$ is not an exact
Abelian sub-category: It does not necessarily contain the kernels, images
and cokernels of its morphisms. However, we have the following result, due
to Ringel:

\begin{prop} \label{p:abcat}
Let $\mfam$ be a family of non-zero, pairwise non-isomorphic left
$R$-modules. Then $\modf \subseteq \rmod$ is a full, extension closed,
exact Abelian subcategory and $\mfam$ is the set of simple objects in
$\modf$, up to isomorphism, if and only if the following conditions hold:
\begin{enumerate}
\item $\enm_R(M_{\alpha})$ is a division ring for all $\alpha \in I$,
\item $\hmm_R(M_{\alpha},M_{\beta}) = 0$ for all $\alpha, \beta \in I$
with $\alpha \ne \beta$.
\end{enumerate}
\end{prop}
\begin{proof}
This follows from Ringel \cite{ri76}, theorem 1.2 and the comments
preceding the theorem.
\end{proof}

Following Ringel, we shall say that an object $M_{\alpha}$ in $\mfam$
is a \emph{point} if $\enm_R(M_{\alpha})$ is a division ring, and that
a pair $(M_{\alpha},M_{\beta})$ of objects in $\mfam$ are \emph{orthogonal}
if $\hmm_R(M_{\alpha},M_{\beta}) = \hmm_R(M_{\beta},M_{\alpha}) = 0$.

If $\mfam$ is a family of orthogonal points, we see that $\modf$ is a
\emph{length category} in the sense of Gabriel \cite{ga73}. In the rest 
of this paper, we shall assume that $\mfam$ is a family of orthogonal
points, unless otherwise specified.

To simplify notation, we shall sometimes say that a left $R$-module $M$
is an iterated extension of $\mfam$ when $M$ is an object of $\modf$.
Since $\modf$ is a length category, there is a Jordan-H\"older theorem
for $\modf$. Consequently, every cofiltration $C$ of $M$ has the same
length $n$ and the same (composition) factors $K_1, \dots, K_n$, up to a
permutation.

We say that an object $M$ in $\modf$ is \emph{uniserial} if the lattice 
of sub-modules of $M$ is a chain (which is the unique composition series 
of $M$). If $M$ is uniserial, then clearly $M$ is indecomposable. But in 
general, the class of indecomposable modules in $\modf$ is larger than 
the class of uniserial modules. We shall later give a characterization 
of when these classes coincide.

\section{Ordered quivers and extension types}
\label{s:extype}

Let $R$ be an associative ring, and let $\mfam = \{ M_{\alpha}: \alpha \in
I \}$ be a family of orthogonal points. We shall define the \emph{extension
type} $\Gm$ of an iterated extension of the family $\mfam$. To give
the extension type is equivalent to giving the length $n$ and the order
vector $\vect \alpha$, and we may therefore consider the extension type as
a (discrete) invariant.

A \emph{quiver} or \emph{directed graph} is a graph $\Gm$, given by a set
$N$ of vertices, a set $E$ of arrows, and maps $s,e: E \to N$. The maps $s,e$
define the starting node $s(a)$ and the ending node $e(a)$ of each arrow
$a \in E$, and we picture $a$ as an arrow from node $s(a)$ to node $e(a)$.
A quiver is said to be \emph{finite} if $N$ and $E$ are finite sets, and
\emph{connected} if the underlying graph is connected. We shall only consider 
quivers which are finite and connected.

An \emph{ordered quiver} is a quiver $\Gm$ together with a total order on
the set $E$ of edges of $\Gm$, such that $e(a)=s(b)$ whenever $a,b \in E$
and $b$ is the successor of $a$. Recall that $b$ is a successor of $a$ if
$a<b$ and the set $\{ c \in E: a<c<b \}$ is empty. To fix notation, we
shall sometimes write $N= \{ 1,2, \dots, p \}$ and $E= \{ a_{12}, a_{23},
\dots, a_{n-1,n} \}$, where $p$ is the number of nodes and $n-1$ is the
number of edges in $\Gm$. The underlying total order of $E$ is given by
$a_{12} < a_{23} < \dots < a_{n-1,n}$, and the definition of an ordered
quiver dictates that $1 \le p \le n$ and that $e(a_{i-i,i})=s(a_{i,i+1})$
for $2 \le i \le n-1$.

Let $(M,C)$ be an object in $\extn$ of length $n$ and with order vector
$\vect \alpha \in I^n$. We let $I(M,C) = \{ \alpha(i): 1 \le i \le n \}$
be the minimal subset $I(M,C) \subseteq I$ such that $(M,C)$ is an
iterated extension of the family $\subfam = \{ M_{\alpha}: \alpha
\in I(M,C) \}$. We define the \emph{extension type} $\Gm$ of the
iterated extension $(M,C)$ to be the ordered quiver given by $N =
I(M,C)$, $E= \{ a_{i-1,i}: 2 \le i \le n \}$, and $s(a_{i-1,i})=
\alpha(i-1)$, $e(a_{i-1,i})=\alpha(i)$.

We remark that the extension type $\Gm$ only depends upon the length $n$
and the order vector $\vect \alpha$ of $(M,C)$, so isomorphic extensions
of extensions of $\mfam$ have the same extension type.

As an example, let us draw all the different extension types of extensions
of extensions of length $n=3$, the first interesting case. When $n=3$, we
must have $1 \le p \le 3$, and there are $5$ isomorphism classes of ordered
directed graphs:

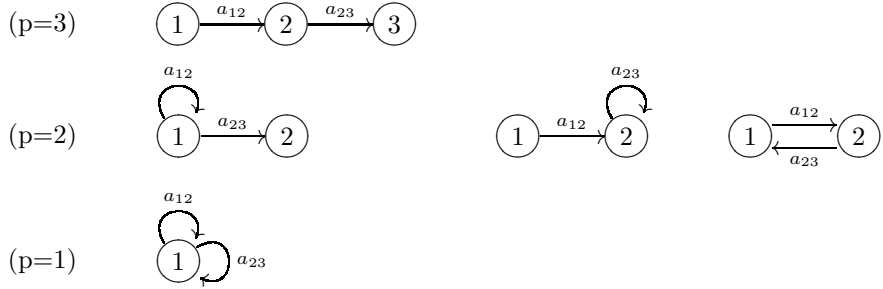
\entrymodifiers={++[o][F-]}
\UseComputerModernTips

\begin{align*}
\tag{p=3}
&\xymatrix{1 \ar[r]^{a_{12}} & 2 \ar[r]^{a_{23}} & 3 } &
&\xymatrix{} &
&\xymatrix{} \\
\tag{p=2}
&\xymatrix{1 \ar@(ul,ur)^{a_{12}} \ar[r]^{a_{23}} & 2 } &
&\xymatrix{1 \ar[r]^{a_{12}} & 2 \ar@(ul,ur)^{a_{23}} } &
&\xymatrix{1 \ar@<1ex>[r]^{a_{12}} & 2 \ar@<1ex>[l]^{a_{23}} } \\
\tag{p=1}
&\xymatrix{1 \ar@(ul,ur)^{a_{12}} \ar@(ur,dr)^{a_{23}} } &
&\xymatrix{} &
&\xymatrix{} \\
\end{align*}

\entrymodifiers={++}

Let $k$ be an algebraically closed field, and let $\Gm$ be an ordered quiver
with $p$ vertices and $n-1$ arrows. There is a $k$-algebra $k[\Gm]$ associated 
with the ordered quiver $\Gm$. This algebra is in a natural way an object of 
the category $\acat p$, the category of complete Artinian $p$-pointed algebras.
We shall briefly recall the definition of $\acat p$.

The category $\Acat p$ is the category of $p$-pointed $k$-algebras: An object
of $\Acat p$ is an associative ring $S$ with structural ring homomorphisms
\[ \xymatrix{k^p \ar[r]^f & S \ar[r]^g & k^p}, \]
such that that composition $g \circ f = \id$, and the morphisms in $\Acat p$
are ring homomorphisms such that the natural diagrams commute (that is, ring
homomorphisms $\phi: S \to S'$ such that $\phi \circ f=f'$ and $g' \circ \phi
= g$). For each object $S$ in $\Acat p$, we denote by $I=I(S)=\ker(g)$ the
\emph{radical ideal} of $S$. The category $\acat p$ is the full sub-category
of $\Acat p$ consisting of objects $S$ such that $S$ is Artinian and complete
in the $I$-adic topology.

We recall some basic facts about $\acat p$: If $S$ is an object of $\Acat p$,
then $S$ is in $\acat p$ if and only if $S$ has finite dimension as vector
space over $k$ and the radical $I$ is nilpotent. In this case, $I$ is the
Jacobson radical of $S$. We denote by $\acat p(n)$ the full sub-category of
$\acat p$ consisting of objects $S$ such that $I(S)^n=0$. Furthermore, any
object $S$ in $\Acat p$ is a matrix ring in the following sense: Denote by
$e_1, \dots, e_p$ the idempotents $e_i = (0, \dots, 1, \dots, 0) \in k^p$ for
$1 \le i \le p$, and let $S_{ij} = e_i S e_j$ for all $1 \le i,j \le p$. Then
$S = \osum S_{ij}$, and $S_{ij} S_{kl} \subseteq \delta_{jk} S_{il}$.
We write $S = (S_{ij})$.

Let $\Gm$ be an ordered quiver with $p$ vertices, and let the edges of $\Gm$ 
be denoted $\{ a_{12}, \dots, a_{n-1,n} \}$ as usual. We define $k[\Gm]$ to be
the object in $\acat p(n)$ given by generators $x_{i-1,i} \in k[\Gm]_{l(i),
l(i-1)}$, with the relations
\begin{equation} \label{e:gmrel}
x_{j-1,j} x_{i-1,i} = 0 \text{ unless } i < j
\end{equation}
for $2 \le i,j \le n$. It follows that $k[\Gm]$ is a finite dimensional vector
space over $k$: It has a natural basis consisting of the non-zero monomials in
$\{ x_{12}, \dots, x_{n-1,n} \}$ of length at most $n-1$ (including $e_1, \dots,
e_p$, which are considered as monomials of length $0$). Furthermore $I^n=0$,
where $I$ is the radical of $k[\Gm]$. So by construction, $k[\Gm]$ is an object
of $\acat p(n)$ when $\Gm$ is an ordered quiver with $p$ vertices and $n-1$ 
arrows.

Let us continue the example of iterated extensions of length $n=3$. In this
case, the algebras $k[\Gm]$ associated to the ordered quivers $\Gm$ shown above
are the following:

\begin{align*}
\tag{p=3}
&\left( \begin{array}{lll}
k e_1 & 0 & 0 \\ k x_{12} & k e_2 & 0 \\ k x_{23} x_{12} & k x_{23} & k e_3 \\
\end{array} \right) & &= &
&\left( \begin{array}{lll}
k & 0 & 0 \\ k & k & 0 \\ k & k & k \\
\end{array} \right) \\
\tag{p=2}
&\left( \begin{array}{ll}
k e_1 + k x_{12} & 0 \\ k x_{23} + k x_{23} x_{12} & k e_2 \\
\end{array} \right) & &= &
&\left( \begin{array}{ll}
k[\epsilon] & 0 \\ k[\epsilon] & k \\
\end{array} \right) \\
&\left( \begin{array}{ll}
k e_1 & 0 \\ k x_{12} + k x_{23} x_{12} & k e_2 + k x_{23} \\
\end{array} \right) & &= &
&\left( \begin{array}{ll}
k & 0 \\ k[\epsilon] & k[\epsilon] \\
\end{array} \right) \\
&\left( \begin{array}{ll}
k e_1 + k x_{23} x_{12} & k x_{23} \\ k x_{12} & k e_2 \\
\end{array} \right) & &= &
&\left( \begin{array}{ll}
k[\epsilon] & (\epsilon) \\ (\epsilon) & k \\
\end{array} \right) \\
\tag{p=1}
&\left( \begin{array}{l}
k + k x_{12} + k x_{23} + k x_{23} x_{12} \\
\end{array} \right) & &= &
&\begin{array}{l}
k\{ x_{12},x_{23} \} / ( x_{12}^2, x_{23}^2, x_{12} x_{23} ) \\
\end{array} \\
\end{align*}

Notice that for each $\Gm$, we have given two different
descriptions of the algebra $k[\Gm]$ in $\acat p$: To the left, we
indicate the natural $k$-linear basis of $k[\Gm]$, and to the
right, we give the multiplicative structure of $k[\Gm]$ (recall
that $\epsilon^2=0$). The multiplicative structure can be worked
out from the natural $k$-linear basis and equation
(\ref{e:gmrel}).

Finally, let us mention that an ordered quiver can be considered as a quiver
with relations: Indeed, let $\Gm$ be the underlying quiver of an ordered
quiver, and consider the relations $a_{j-1,j} a_{i-1,i} = 0$ for all $i,j$
with $2 \le j \le i \le n$.

\section{Noncommutative deformations of modules}
\label{s:defm}

Let $k$ be an algebraically closed field, let $R$ be an associative 
$k$-algebra, and let $\mfam = \{ M_1, \dots, M_p \}$ be a finite 
family of left $R$-modules. There is a deformation functor 
        \[ \defm: \acat p \to \sets, \]
describing the simultaneous deformations of the family $\mfam$ of left
$R$-modules, and a theory of noncommutative deformations of modules 
related to this functor: This theory is due to Laudal, and it is 
described in several preprints, see Laudal \cite{lau95,lau98,lau00}. 
However, we find it more convenient to give references to Eriksen 
\cite{er02}, which is a version of the theory adapted to the study of 
left modules.

Let us briefly recall the definition of the deformation functor $\defm$:
Let $S$ be an object in $\acat p$. A \emph{lifting of the family $\mfam$ 
to $S$} is a left $R \otimes_k \op S$-module $M_S$ together with 
isomorphisms $\eta_i: M_S \otimes_S k_i \to M_i$ of left $R$-modules for 
$1 \le i \le p$, such that $M_S \cong (M_i \otimes_k S_{ij})$ considered 
as right $S$-modules. Recall that $k_i$ is the image of the $i$'th 
projection of $k^p$, which is an $S$-module via $g: S \to k^p$. We let 
$\defm(S)$ denote the set of all equivalence classes of liftings of the 
family $\mfam$ to $S$, where $M_S$ and $M'_S$ are equivalent liftings if 
there is an isomorphism $\tau: M_S \to M'_S$ of left $R \otimes_k \op 
S$-modules such that the natural diagrams commute (that is, such that 
$\eta'_i \circ (\tau \otimes_S k_i) = \eta_i$ for $1 \le i \le p$). 

For the rest of this sections, we shall assume that $\dim_k 
\ext^m_R(M_i,M_j)$ is finite for $1 \le i,j \le p, \; m = 1,2$. In this 
case, there exists a \emph{pro-representable hull} $H = H(\mfam)$ for 
the deformation functor $\defm$ (which is unique up to non-canonical 
isomorphism):

\begin{thm}[Laudal]
Let $k$ be an algebraically closed field, let $R$ be an associative 
$k$-algebra, and let $\mfam = \{ M_1, \dots, M_p \}$ be a finite family of 
left $R$-modules such that $\dim_k \ext^m_R(M_i,M_j)$ is finite for $1 \le 
i,j \le p, \; m = 1,2$. Then there exists a pro-representable hull $H$ for 
the deformation functor $\defm: \acat p \to \sets$.
\end{thm}
\begin{proof}
The hull $H$ can be constructed along well-known lines, via the obstruction 
morphism $o: \obstrdef \to \infdef$ (see Eriksen \cite{er02}, theorem 4.2),
or via (non-symmetric) matric Massey products (outlined in Laudal 
\cite{lau00}).  
\end{proof}

Recall that a pro-representable hull $H$ for the functor $\defm$ is an object
of the pro-category $\cacat p$ such that there exists a smooth morphism of
functors on $\acat p$
        \[ \mor(H,-) \to \defm, \] 
which is an isomorphism when restricted to a morphism of functors on 
$\acat p(2)$. The pro-category $\cacat p$ is the full sub-category of $\Acat p$ 
consisting of objects $S$ which are complete in the $I$-adic topology and such 
that $S_n = S / I^n$ is an object in $\acat p(n)$ for all $n \ge 1$. 

\begin{prop} \label{p:affstr}
Let $H$ be the pro-representable hull of the deformation functor $\defm: \acat p 
\to \sets$, and let $S$ be any object in $\acat p$. Then $\mor(H,S)$ has a 
natural structure as an affine scheme over $k$. In particular, $\mor(H,k[\Gm])$ 
is an affine scheme over $k$ for any ordered quiver $\Gm$ with $p$ nodes. 
\end{prop}
\begin{proof}
Let $V_{ij}^m = \ext^m_R(M_j,M_i)^*$ for $1 \le i,j \le p, \; m = 1,2$, and 
choose $k$-linear bases $\{ s_{ij}(\alpha): 1 \le \alpha \le d_{ij} \}$ for 
$V_{ij}^1$ and $\{ t_{ij}(\beta): 1 \le \beta \le r_{ij} \}$ for $V_{ij}^2$. 
Then $\infdef$ is the (free) formal matrix ring with generators 
$\{ s_{ij}(\alpha) \}$, $\obstrdef$ is the (free) formal matrix ring with 
generators $\{ t_{ij}(\beta) \}$, and the obstruction morphism is a morphism 
$o: \obstrdef \to \infdef$ in $\cacat p$. Let $f_{ij}(\beta) = o(t_{ij}(\beta)) 
\in \infdef_{ij}$ for $1 \le i,j \le p, \; 1 \le \beta \le r_{ij}$. There is 
a natural surjection $\infdef \to H$, and its kernel is generated by the 
relations $\{ f_{ij}(\beta): 1 \le i,j \le p, \; 1 \le \beta \le r_{ij} \}$. 
Clearly this surjection induces an injective map of sets $\mor(H,S) \to 
\mor(\infdef,S)$. Since $I(S)^n=0$ for some $n \ge 1$, we have that 
        \[ \mor(\infdef,S) = \mor(\infdef_n,S) = \prod_{i,j} 
        \hmm_k(V_{ij}^1,W_{ij}), \]
where $W_{ij} = I(S)_{ij}$ with basis $\{ w_{ij}(\gamma): 1 \le \gamma \le 
v_{ij} \}$. So $\mor(\infdef,S) \cong \aff^N$, where $N = \sum d_{ij} v_{ij}$. 
We obtain a set of coordinates $\{ z_{ij}(\alpha,\gamma) \}$ for $\aff^N$, 
where the coordinate $z_{ij}(\alpha,\gamma)$ corresponds to a morphism 
$\phi_{ij}(\alpha,\gamma) \in \mor(\infdef,S)$ given by 
        \[ \phi_{ij}(\alpha,\gamma)(s_{ij}(\alpha')) = \delta_{\alpha,\alpha'} 
        w_{ij}(\gamma) \]
for $1 \le i,j \le p, \; 1 \le \gamma \le v_{ij}, \; 1 \le \alpha,\alpha' \le 
d_{ij}$. Let $\phi=(a_{ij}(\alpha,\gamma)) \in \mor(\infdef,S)$, then $\phi 
\in \mor(H,S)$ if and only if $\phi(f_{ij}(\beta))=0$ for all $1 \le i,j \le p, 
\; 1 \le \beta \le r_{ij}$. But we have 
        \[ \phi(f_{ij}(\beta)) = \sum_{\gamma} 
        f_{ij}(\beta)(\{ a_{ij}(\alpha,\gamma): 1 \le \alpha \le d_{ij} \}) \;  
        w_{ij}(\gamma), \]
so $\phi(f_{ij}(\beta))=0$ if and only if we have the equations
        \[ f_{ij}(\beta)(\{ a_{ij}(\alpha,\gamma): 1 \le \alpha \le d_{ij} \}) 
        = 0 \]
for $1 \le i,j \le p, \; 1 \le \beta \le r_{ij}, \; 1 \le \gamma \le v_{ij}$. 
Notice that $\phi(I(\infdef)^n)=0$, so the above equations corresponds to 
polynomial relations $R_{ij}(\beta,\gamma) \in k[ \{ z_{ij}(\alpha,\gamma) \}]$, 
and therefore $\mor(H,S)$ is the affine sub-scheme of $\mor(\infdef,S) \cong 
\aff^N$ defined by those relations. 
\end{proof}
 
We denote by $\x$ the affine scheme $\mor(H,k[\Gm])$. Moreover, we denote by 
$M_H$ be the \emph{versal family} $M_H \in \defm(H)$ corresponding to the smooth 
morphism $\mor(H,-) \to \defm$ via Yonedas lemma. For each point $\phi \in \x$, 
there is a deformation $M_{\phi} = \defm(\phi)(M_H) \in \defm(k[\Gm])$. By 
construction, the morphism $\mor(H,-) \to \defm$ is smooth, so the map of sets  
        \[ \x \to \defm(k[\Gm]), \]
given by $\phi \mapsto M_{\phi}$, is surjective.  

We shall explain how to calculate the surjection above in concrete terms: Let 
$\phi \in \x$, so $\phi: H \to k[\Gm]$ is a morphism in $\cacat p$. Let $M_H 
\in \defm(H)$ be the versal family defined over $H$. Then the deformation 
$M_{\phi} \in \defm(k[\Gm])$ is given by 
        \[ M_{\phi} = (M_i \otimes_k k[\Gm]_{ij}) \]
considered as a right $k[\Gm]$-module, and the left $R$-module structure of 
$M_{\phi}$ is determined by
\begin{equation} \label{e:phidef}
r (m_i \otimes e_i) = (\id \otimes \phi) (r(m_i \otimes e_i)) 
\end{equation} 
for all $r \in R, \; m_i \in M_i$. The expression $r (m_i \otimes e_i)$ on the 
right hand side of equation (\ref{e:phidef}) is the left multiplication of $r 
\in R$ with $m_i \otimes e_i$ considered as an element of $M_H$. This makes
it possible to compute $M_{\phi} \in \defm(k[\Gm])$ when $\phi \in \x$ is 
given, assuming that the versal family $M_H$ can be computed. 

For the rest of this section, assume that $\mfam = \{ M_1, \dots, M_p \}$ is 
a finite family of non-zero, pairwise non-isomorphic left $R$-modules, and 
consider the category $\extn$ of iterated extensions of the family 
$\mfam$. For any ordered quiver $\Gm$ with vertices $N = \{ 1,2, \dots, p \}$, 
we denote by $\extncl(\mfam,\Gm)$ the set of isomorphism classes of extensions 
of extensions of the family $\mfam$ with extension type $\Gm$. It is clear 
that the forgetful functor $\extn \to \modf$ maps $\extncl(\mfam,\Gm)$ to a 
set of isomorphism classes of left $R$-modules, which we shall denote by 
$\modfcl(\mfam,\Gm)$. Moreover, the above construction defines a natural 
surjective map $\extncl(\mfam,\Gm) \to \modfcl(\mfam,\Gm)$. 

\begin{thm}[Laudal] \label{t:defext}
Let $\mfam = \{ M_1, \dots, M_p \}$ be a finite family of non-zero, pairwise
non-isomorphic left $R$-modules such that $\dim_k \ext^m_R(M_i,M_j)$ is finite 
for $1 \le i,j \le p, \; m=1,2$, and let $\Gm$ be an ordered quiver with 
vertices $N = \{ 1,2, \dots, p \}$. Then there is a natural bijection between 
$\defm(k[\Gm])$ and $\extncl(\mfam,\Gm)$. 
\end{thm}
\begin{proof}
Let the $R$-$k[\Gm]$ bimodule $(M_i \otimes_k k[\Gm]_{ij})$ be a lifting of the 
family $\mfam$ to $k[\Gm]$. We show how to construct an iterated extension
$(M,C)$ of the family $\mfam$ with extension type $\Gm$: We let $M' = (M_i 
\otimes_k k[\Gm]_{i,\alpha(1)}) \subseteq (M_i \otimes_k k[\Gm]_{ij})$, this is 
by construction invariant under left multiplication by $R$. Consider sequences 
$j = (j_1, \dots, j_r)$ with $1 \le j_1 < j_2 < \dots < j_r \le n-1$ and 
$\alpha(j_h+1)=\alpha(j_{h+1})$ for $1 \le h \le r-1$. We say that $j$ is 
broken if $j_h+1 \not = j_{h+1}$ for some $h$, otherwise $j$ is unbroken. Let us 
denote by $M'(j)$ the $k$-linear space $M_{j_r+1} \otimes x_{j_r,j_r+1} \dots 
x_{j_1,j_1+1}$, by $M'(B)$ the sum of all $k$-linear spaces $M'(j)$ with $j$ 
broken, and $M(U)$ the sum of all $k$-linear spaces $M'(j)$ with $j$ unbroken. 
Notice that $M'(B)$ is invariant under left multiplication with $R$. We define 
$M = M' / M'(B)$, which has a natural left $R$-module structure, and clearly $M 
\cong M(U)$ considered as a $k$-linear space. For $0 \le i \le n$, let $F_i$ be 
the sum of all $k$-linear spaces $M'(j)$ where $j$ is unbroken of length at 
least $i$. Then $F_i \subseteq M$ are also invariant under left multiplication 
with $R$, so $C$ is a co-filtration of $M$ when $C_i = M / F_i$ for $0 \le i \le 
n$. Clearly, this co-filtration satisfies $K_i = \ker(C_i \to C_{i-1}) \cong 
M_{\alpha(i)}$ for $1 \le i \le n$. It follows that $(M,C)$ is an extension of 
extensions of the family $\mfam$ of extension type $\Gm$. It is also 
straight-forward to check that equivalent liftings of $\mfam$ to $k[\Gm]$ gives 
isomorphic iterated extensions: Any isomorphism between liftings will map 
$M'$ to $M'$, $M'(B)$ to $M'(B)$, and $F_i$ to $F_i$. So we have constructed a 
well-defined map from $\defm(k[\Gm])$ to $\extncl(\mfam,\Gm)$.

It only remains to see that this map is a bijection: But given an extension of 
extensions $(M,C)$ of the family $\mfam$, consider $(M_i \otimes_k k[\Gm]_{ij})$ 
as a right $k[\Gm]$-module. It is easy to see that the left $R$-module structure 
of the co-filtration $C$ on $M$ will generate a left $R$-module structure on 
$(M_i \otimes_k k[\Gm]_{ij})$ compatible with the right $k[\Gm]$-module 
structure. Furthermore, isomorphic iterated extensions give equivalent 
liftings of $\mfam$ to $k[\Gm]$. So we have constructed an inverse to the map 
of sets described above. 
\end{proof}

\begin{cor}
There are natural surjections $\x \to \extncl(\mfam,\Gm) \to \modfcl(\mfam,\Gm)$. 
In particular, the sets $\extncl(\mfam,\Gm)$ and $\modfcl(\mfam,\Gm)$ are 
quotients of the affine scheme $\x$.
\end{cor}

We remark that these quotients are computable, in principle: In proposition 
\ref{p:affstr}, we have shown how to construct the affine scheme $\x$ when 
the hull $H$ of $\defm$ is known. Moreover, we have shown above how to 
calculate the surjection $\x \to \defm(K[\Gm])$ when the versal family $M_H 
\in \defm(H)$ is known. The identification $\defm(k[\Gm]) \cong 
\extncl(\mfam,\Gm)$ is explicitly given in the proof of theorem \ref{t:defext}, 
and the surjection $\extncl(\mfam,\Gm) \to \modfcl(\mfam,\Gm)$ is natural, 
induced by the forgetful functor $(M,C) \mapsto M$.

\section{Species}

We say that $S = (K_{\alpha}, E_{\alpha, \beta})$ is a \emph{species} 
indexed by the set $I$ if $K_{\alpha}$ is a division ring and $E_{\alpha,
\beta}$ is a $K_{\beta}- K_{\alpha}$ bimodule for all $\alpha, \beta \in 
I$. Let $k$ be a fixed commutative field. We say that $S$ is a $k$-species
if $k$ is contained in $K_{\alpha}$ for all $\alpha \in I$ in such a way 
that $c \xi = \xi c$ for all $c \in k, \; \xi \in E_{\alpha,\beta}$ and 
$\dim_k K_{\alpha}$ is finite. Moreover, $S$ is called a $k$-quiver if
in addition $K_{\alpha} = k$ for all $\alpha \in I$. 

If $S$ is a $k$-quiver, it is completely determined by the \emph{Gabriel
quiver} of $S$, formed in the following way: The set of vertices of the
Gabriel quiver is $I$, and for each pair of vertices $\alpha, \beta \in I$, 
there are $\dim_k E_{\alpha, \beta}$ arrows from $\alpha$ to $\beta$.
 
Let $\mfam$ be a family of orthogonal points in $\rmod$ indexed by $I$, and 
consider the corresponding length category $\modf$. We define the species of 
$\modf$ to be the species indexed by $I$ given by $K_{\alpha} = 
\enm_R(M_{\alpha})$ and $E_{\alpha,\beta} = \ext^1_R(M_{\alpha},M_{\beta})$ 
for all $\alpha, \beta \in I$. 

In what follows, we shall be particularly interested in the case when $k$
is an algebraically closed commutative field, $R$ is a $k$-algebra, and 
$\mfam$ is a family of orthogonal points in $\rmod$ such that 
$\enm_R(M_{\alpha}) = k$ for all $\alpha \in I$. In this case, the species 
of $\modf$ is clearly a $k$-quiver. 

It is well-known that the species of the length category $\modf$ contains a 
lot of information about the category, see for instance Gabriel \cite{ga73}.
This fact can be explained by a general principle: The length category 
$\modf$ is completely determined by its species $S$ and the obstruction 
theory of the family $\mfam$ of simple objects. In theorem \ref{t:defext},
we have proved this principle under some finiteness conditions, using 
non-commutative deformations of modules.

\section{The hereditary case}

In the previous section, we have seen that the length category
$\modf$ is determined by its species $S$ and the obstruction
theory of $\mfam$. We shall consider the \emph{unobstructed}
cases, which are the easiest ones. That is, we shall consider
the cases in which $\modf$ is completely determined by its 
species $S$. 

We say that $\modf$ is \emph{hereditary} if $\ext^2_R(M_{\alpha},
M_{\beta})=0$ for all $\alpha, \beta \in I$. This is clearly
a sufficient (but not necessary) condition for the length 
category $\modf$ to be unobstructed. To understand how $\modf$
is related to the species $S$ in this case, we shall use some
results of Deng, Xiao \cite{de-xi98}:

In the rest of this section, let $k$ be an algebraically closed 
commutative field, let $R$ be a $k$-algebra, and let $\mfam$ be a 
family of orthogonal points in $\rmod$ such that $\modf$ is a 
hereditary category and $\enm_R(M_{\alpha})=k$ for all $\alpha 
\in I$. We shall also assume that the following finiteness 
conditions hold:

\begin{itemize}
\item $I$ is a finite set,
\item $\dim_k \ext^1_R(M_{\alpha},M_{\beta})$ is finite for all 
$\alpha, \beta \in I$.
\end{itemize}

The species $S$ of the hereditary length category $\modf$ is clearly 
a $k$-quiver, so it is completely determined by the corresponding 
Gabriel quiver $Q$. The finiteness conditions above means that
the Gabriel quiver is finite.

We consider the \emph{finite representations} $V$ of $Q$: These 
consist of a finite dimensional vector space $V_{\alpha}$ over $k$
for all $\alpha \in I$, and a $k$-linear map $V_{a}: V_{\alpha} \to
V_{\beta}$ for all arrows $a: \alpha \to \beta$. We say that a finte
representation $V$ is \emph{small} or \emph{nilpotent} if there is a 
positive integer $n \ge 1$ such that $V_p = 0$ for all paths $p$ of 
length $n$ in $Q$.  

\begin{thm}[Deng-Xiao]
Let $\mfam$ be a family of left $R$-modules satisfying the conditions
above. Then there is a natural exact equivalence of categories between 
$\modf$ and the category of small representations of Gabriel quiver $Q$. 
\end{thm}
\begin{proof}
See Deng, Xiao \cite{de-xi98}, theorem 1.5.
\end{proof}

In other words, the category $\modf$ is completely determined by the 
Gabriel quiver $Q$, and therefore also by the species of $\modf$. Notice
that if $Q$ is loop-free, then all finite representations are small. In
this case, we can tell when the category $\modf$ is finite, tame or wild
by considering the well-known classification of quivers into these 
classes. 

An interesting question is the following: Let $Q$ be a finite quiver
with loops. When is the category of small representations of $Q$ finite,
tame and wild? We do not know if the answer to this question is known.

\section{Indecomposable and uniserial objects}

Let $k$ be an algebraically closed field, let $R$ be an associative
$k$-algebra, and let $\mfam = \{ M_{\alpha}: \alpha \in I \}$ be a
family of orthogonal points. We are interested in the indecomposable 
modules in the length category $\modf$.

Since $\modf$ is a length category, we have a
Krull-Remak-Schmidt-Azumaya theorem. So every object in $\modf$ has 
a finite indecomposable decomposition, unique up to a permutation, 
and $\enm_R(M)$ is a local ring for all indecomposable objects $M 
\in \modf$. We say that $\modf$ is a \emph{uniserial category} if 
every indecomposable object in $\modf$ is uniserial.

\begin{lem} \label{l:induni}
Let $M$ be an object of $\modf$, and consider the following conditions:
\begin{enumerate}
\item $M$ is uniserial,
\item $M$ has a unique minimal submodule,
\item $M$ is indecomposable.
\end{enumerate}
Then we have $1) \Rightarrow 2) \Rightarrow 3)$. In particular, all
conditions are equivalent if and only if $\modf$ is a uniserial category.
\end{lem}
\begin{proof}
The implication $1) \Rightarrow 2)$ is obvious. If $M = N_1 \oplus N_2$
is a direct decomposition with $N_i \ne 0$ for $i=1,2$, there are minimal
sub-modules $K_i \subseteq N_i$ for $i=1,2$. This shows that $2)
\Rightarrow 3)$. The last part is clear.
\end{proof}

The implication $3) \Rightarrow 1)$ in the above lemma clearly holds if
$M$ has length $n=2$. But already in the case $n=3$, it is very easy
to come up with counterexamples:

\begin{lem} \label{l:ce1}
Let $\mfam$ be a family of orthogonal points. If $\mfam$ contains modules
$S,T$ such that $\enm_R(K)=k$ for $K=S,T$ and $\dim_k \ext^1_R(S,T)
\ge 2$, then there exists an object $M$ in $\modf$ of length $n=3$ such
that $M$ is a non-uniserial $R$-module with a unique simple submodule.
\end{lem}
\begin{proof}
From McConnel, Robson \cite{mc-ro73}, proposition 3.3 it follows that there
exists non-split extensions $U,V$ of $S$ by $T$ such that $U$ and $V$ are
not isomorphic as $R$-modules. Let $M$ be the cokernel of the diagonal map
$T \to U \oplus V$, then $M$ is a non-uniserial module with the unique
simple submodule $T$ by McConnell, Robson \cite{mc-ro73}, proposition 6.1.
\end{proof}

\begin{lem} \label{l:ce2}
Let $\mfam$ be a family of orthogonal points. If $\mfam$ contains modules
$S,T,U$ such that $\enm_R(K)=k$ for $K=S,T,U$, $\ext^1_R(U,S),
\ext^1_R(U,T) \ne 0$, and $S,T$ are non-isomorphic, then there exists an
object $M$ in $\modf$ of length $n=3$ such that $M$ is indecomposable and
such that $S$ and $T$ are simple submodules of $M$.
\end{lem}
\begin{proof}
Let $\xi_1, \xi_2 \ne 0$ be extensions of $U$ by $S$ and $U$ by $T$, and let
$\psi_i$ be representatives of $\xi_i$ in Hochschild cohomology for $i=1,2$.
Let furthermore $M$ be the extension of $U$ by $S \oplus T$ given by $(\xi_1,
\xi_2)$. Then $M \cong S \oplus T \oplus U$ as a vector space over $k$, and
the $R$-module structure of $M$ is defined by the representatives $\psi_1,
\psi_2$. A calculation shows that $\enm_R(M)=k$ if $U$ is not isomorphic to
any of $S,T$, and that $\enm_R(M)=k[x]/(x^2)$ if $U$ is isomorphic to one of
the modules $S,T$. In either case, $\enm_R(M)$ is a local ring, and therefore
$M$ is indecomposable. Since $S \oplus T \subseteq M$, it is clear that $M$
has simple submodules $S,T$.
\end{proof}

\begin{lem} \label{l:ce3}
Let $\mfam$ be a family of orthogonal points. If $\mfam$ contains modules
$S,T,U$ such that $\enm_R(K)=k$ for $K=S,T,U$, $\ext^1_R(S,U),
\ext^1_R(T,U) \ne 0$, and $S,T$ are non-isomorphic, then there exists an
object $M$ in $\modf$ of length $n=3$ such that $M$ is indecomposable but
not uniserial.
\end{lem}
\begin{proof}
The proof is similar to the proof of the lemma \ref{l:ce2}: We consider
$M = U \oplus S \oplus T$ considered as a vector space over $k$, and let the
$R$-module structure of $M$ be given by non-split extension $\xi_S \in
\ext^1_R(S,U)$ and $\xi_T \in \ext^1_R(T,U)$ via Hochschild cohomology.
A calculation shows that $\enm_R(M)=k$ if $U$ is not isomorphic to $S$ or
$T$, and $\enm_R(M) \cong k[x]/(x^2)$ otherwise, so $M$ is indecomposable
in both cases. On the other hand, $M$ has to submodules of length $2$, so
$M$ is not uniserial.
\end{proof}

\begin{cor}
Let $\mfam$ be a family of orthogonal points such that the species of
$\mfam$ is a $k$-quiver. If $\modf$ is a uniserial category, then we
have
\begin{enumerate}
\item $\sum_{\beta \in I} \dim_k (\ext^1_R(M_{\alpha},M_{\beta})) \le 1$
for all $\alpha \in I$,
\item $\sum_{\alpha \in I} \dim_k (\ext^1_R(M_{\alpha},M_{\beta})) \le 1$
for all $\beta \in I$.
\end{enumerate}
\end{cor}
\begin{proof}
This follows from lemma \ref{l:ce1}, \ref{l:ce2} and \ref{l:ce3}.
\end{proof}

We shall later see that these conditions are also sufficient for $\modf$
to be a uniserial category. We remark that this criterion has been known
since the 60's, see Amdal, Ringdal \cite{am-ri68} and Gabriel \cite{ga73},
section 8.3.

\section{The uniserial case}
\label{s:uniserial}

Let $k$ be an algebraically closed commutative field, let $R$ be an
associative $k$-algebra, and let $\mfam$ be a family of orthogonal 
points in $\rmod$ such that the species of $\mfam$ is a $k$-quiver.
We shall determine when $\modf$ is an uniserial category. In the 
process, we shall also show that if $\modf$ is uniserial then it is 
unobstructed, regardless if it is hereditary or not.

We recall that the species of $\mfam$ is a $k$-quiver if and only if
$\enm_R(M_{\alpha}) \cong k$ for all $\alpha \in I$, which is
equivalent to the condition that for any $\alpha \in I$ and any
endomorphism $\phi \in \enm_R(M_{\alpha})$, $\phi$ is algebraic over
$k$. By Quillen's lemma, this is the case when $M$ is any simple
module over a ring $R=\diff(A)$ of $k$-linear differential operators
on $A$, such that $\gr \diff(A)$ is a finitely generated $k$-algebra,
see Quillen \cite{qu69}. In particular, any family $\mfam$ of simple
modules over the first Weyl algebra $R=\weyl$ satisfy these 
conditions.

In the first part of this section, we shall assume that the family
$\mfam$ satisfy the following additional conditions:
\begin{itemize}
\item[(*)] $\sum_{\beta \in I} \dim_k (\ext^1_R(M_{\alpha},M_{\beta}))
\le 1$ for all $\alpha \in I$,
\item[] $\sum_{\alpha \in I} \dim_k (\ext^1_R(M_{\alpha},M_{\beta}))
\le 1$ for all $\beta \in I$.
\end{itemize}
When $\mfam$ satisfy (*), we shall classify all iterated extensions
$(M,C)$ in $\extn$ such that $M$ is indecomposable, up to isomorphism in
$\extn$, and all indecomposable modules $M \in \modf$ up to isomorphism.
However, it is useful to start looking at iterated extensions $(M,C)$
such that $\xi_i \ne 0$ for $2 \le i \le n$:

\begin{lem} \label{l:indextnstr}
Let $\mfam$ be a family satisfying (*), and let $(M,C)$ be an iterated
extension of the family $\mfam$ of length $n \ge 2$ such that $\xi_i
\ne 0$ for $2 \le i \le n$. For any module $K \in \mfam$, the map
$\ext^1_R(C_{i-1},K) \to \ext^1_R(K_{i-1},K)$ induced by the inclusion
$K_{i-1} \subseteq C_{i-1}$ is an isomorphism for $2 \le i \le n$.
\end{lem}
\begin{proof}
We show the result by induction on $n$. Since $C_1=K_1$ by definition,
the result is clearly true for $n=2$. So let $n \ge 3$, and assume that
the result holds for all integers less than $n$. In particular, this
implies that
\[ \ext^1_R(C_{n-2},K_{n-1}) \to \ext^1_R(K_{n-2},K_{n-1}) \]
is an isomorphism, and consequently $\ext^1_R(K_{n-2},K_{n-1}) \ne 0$.
Consider the long exact sequence of the functor $\hmm_R(-,K)$ applied to
the extension $\xi_{n-1}$. Since $\xi_{n-1} \ne 0$ and $\enm_R(K_{n-1})
\cong k$, it follows that $\hmm_R(K_{n-1},K) \to \ext^1_R(C_{n-2},K)$ is
injective, and therefore the sequence
\begin{align*}
0 &\to \hmm_R(K_{n-1},K) \to \ext^1_R(C_{n-2},K) \\
&\to \ext^1_R(C_{n-1},K) \to \ext^1_R(K_{n-1},K)
\end{align*}
is exact. But $\hmm_R(K_{n-1},K) \cong \ext^1_R(K_{n-2},K)$, since
$\ext^1_R(K_{n-2},K) \ne 0$ if and only if $K \cong K_{n-1}$. So
$\hmm_R(K_{n-1},K) \to \ext^1_R(C_{n-2},K)$ is an isomorphism by the
induction hypothesis. So the map $\ext^1_R(C_{n-1},K) \to
\ext^1_R(K_{n-1},K)$ is injective. If $K=K_n$, then this map is also
surjective, since it maps $\xi_n$ to $\tau_n$, and $\xi_n \ne 0$. If
$K \not \cong K_n$, then $\ext^1_R(K_{n-1},K) = 0$ and the map is an
isomorphism as well.
\end{proof}

\begin{cor} \label{c:gdata}
Let $\mfam$ be a family satisfying (*), and let $(M,C)$ be an iterated
extension of the family $\mfam$ such that $\xi_i \ne 0$ for $2 \le i
\le n$. Then $\tau_i \ne 0$ for $2 \le i \le n$. In particular, the
extension type of $(M,C)$ is uniquely determined by $(\alpha(1),n) \in
I \times \n$.
\end{cor}

It is also useful to notice that if $(M,C)$ is an iterated extension
of the type considered above, then $M$ is an indecomposable left
$R$-module:

\begin{lem} \label{l:ind}
Let $\mfam$ be a family satisfying (*), and let $(M,C)$ be an iterated 
extension of $\mfam$ of length $n$ such that $\xi_i \ne 0$ for $2 \le
i \le n$. Then $M$ is indecomposable.
\end{lem}
\begin{proof}
Because of lemma \ref{l:induni}, it is enough to show that $K_n \subseteq
M$ is the unique minimal submodule of $M$: Assume that $K \in \mfam$ and
that $\phi: K \to M$ is injective. Because of lemma \ref{l:indextnstr},
we may assume that $M \cong K_n \oplus \dots \oplus K_1$ as a vector space
over $k$. Moreover, we may assume that the left $R$-module structure of
$M$ is given by
\[ r(k_n, \dots, k_1) = (rk_n + \psi^n_r(k_{n-1}), rk_{n-1} +
\psi^{n-1}_r(k_{n-2}), \dots, rk_2 + \psi^2_r(k_1), rk_1) \]
for all $r \in R, \; (k_n, \dots, k_1) \in M$, where $\psi^i \in
\der_k(R,\hmm_k(K_{i-1},K_i))$ represents $\tau_i \in \ext^1_R(K_{i-1},
K_i)$ in Hochschild cohomology. Clearly, there are $k$-linear maps
$\phi_i: K \to K_i$ such
that $\phi(k)=(\phi_n(k), \dots, \phi_1(k))$, and $\phi_n:K \to K_n$ is
$R$-linear as well. But a simple calculation, using the fact that $\phi$
is $R$-linear and $\enm_R(K_i) \cong k$, shows that $\phi_n \ne 0$. So
$\phi_n(K) = K_n \subseteq M$, and therefore $K_n \subseteq M$ is the
unique minimal submodule of $M$.
\end{proof}

Let $(\alpha,n) \in I \times \n$. We consider the set of order vectors
$\vect \alpha \in I^n$ such that $\alpha(1) = \alpha$ and
$\ext^1_R(M_{\alpha(i-1)},M_{\alpha(i)}) \ne 0$ for $2 \le i \le n$. We
say that the couple $(\alpha,n) \in I \times \n$ is \emph{admissible}
if such an order vector exists. If this is the case, we know from (*)
that the order vector $\vect \alpha$ is uniquely determined by $(\alpha,
n)$, and we shall say that $\vect \alpha$ the order vector associated with
the admissible couple $(\alpha,n)$.

Let $(\alpha,n) \in I \times \n$ be an admissible couple, and let $\vect
\alpha$ be the associated order vector. We denote by $\Gm(\alpha,n)$ the
ordered quiver defined by the order vector $\vect \alpha$. If $\mfam$ is
a family satisfying (*) and $(M,C)$ is an iterated extension of the
family $\mfam$ such that $\xi_i \ne 0$ for $2 \le i \le n$, then the
extension type of $(M,C)$ is $\Gm(\alpha,n)$ with $\alpha = \alpha(1)$ by
corollary \ref{c:gdata}.

Let $(M,C)$ be an iterated extension of the family $\mfam$ of
extension type $\Gm$. To simplify notation, we shall write $(M,C)$ for
the isomorphism class of $(M,C)$ in $\extncl(\mfam,\Gm)$, and $M$ for
the isomorphism class of $M$ in $\modfcl(\mfam,\Gm)$.

We shall denote by $\indextncl(\mfam,\Gm) \subseteq \extncl(\mfam,\Gm)$
the subset of isomorphism classes $(M,C)$ in $\extncl(\mfam,\Gm)$ such
that $M$ is indecomposable, and by $\nsextncl(\mfam,\Gm) \subseteq
\extncl(\mfam,\Gm)$ the subset of isomorphism classes $(M,C)$ such that
$\xi_i \ne 0$ for $2 \le i \le n$. We denote by $\indmodfcl(\mfam,\Gm)$
and $\nsmodfcl(\mfam,\Gm)$ the images of $\indextncl(\mfam,\Gm)$ and
$\nsextncl(\mfam,\Gm)$ under the natural surjection $\extncl(\mfam,\Gm)
\to \modfcl(\mfam,\Gm)$.

Let $(\alpha,n)$ be an admissible couple, and consider the ordered
quiver $\Gm(\alpha,n)$. We shall find the classification spaces
$\nsextncl(\mfam,\Gm)$ and $\nsmodfcl(\mfam,\Gm)$ when $\mfam$ is
family satisfying (*) and $\Gm = \Gm(\alpha,n)$. For this, we only
need the following simple lemma:

\begin{lem} \label{l:extmod}
Let $M,N$ be non-zero left $R$-modules, let $E, E'$ be extensions
of $M$ by $N$, and let $\xi, \xi' \in \ext^1_R(M,N)$ be the
corresponding classes. If $\xi' = \psi \xi \phi$ for automorphisms
$\phi \in \aut_R(M), \; \psi \in \aut_R(N)$, then $E \cong E'$
considered as left $R$-modules. In particular, $E \cong E'$ as
$R$-modules when $\xi' = \alpha \xi$ with $\alpha \in k^*$.
\end{lem}
\begin{proof}
Clearly, we have $k \cong k \id_M \subseteq \enm_R(M)$ and $k^*
\subseteq \aut_R(M)$ when $M$ is a non-zero $R$-module. The rest
follows from McConnel, Robson \cite{mc-ro73}, proposition 3.3 and
the preceding paragraph.
\end{proof}

\begin{prop} \label{p:indscl}
Let $\mfam$ be a family satisfying (*), let $(\alpha,n) \in I \times \n$
be an admissible couple, and let $\Gm = \Gm(\alpha,n)$. Then we have
$\nsextncl(\mfam,\Gm) \cong (k^*)^{n-1}$ and $\nsmodfcl(\mfam,\Gm) = \{*\}$.
\end{prop}
\begin{proof}
Let $\Gm = \Gm(\alpha,n)$, let $K_i = M_{\alpha(i)}$ and choose a basis for
$\ext^1_R(K_{i-1},K_i)$ for $2 \le i \le n$. There is a map
$\nsextncl(\mfam,\Gm) \to (k^*)^{n-1}$ given by the composition $(M,C)
\mapsto (\xi_2, \dots, \xi_n) \mapsto (\tau_2, \dots, \tau_n) \mapsto
(k^*)^{n-1}$. This map is clearly injective, since $\ext^1_R(C_{i-1},K_i)
\to \ext^1_R(K_{i-i},K_i)$ is an isomorphism and the extensions $\xi_i$
determine $(M,C)$. Furthermore, the map is surjective by lemma \ref{l:ind}.
Finally, lemma \ref{l:extmod} shows if $(M,C)$ and $(M',C')$ are any
isomorphism classes in $\nsextncl(\mfam,\Gm)$, then $M \cong M'$ considered
as left $R$-modules.
\end{proof}

Let $(M,C)$ be an iterated extension of the family $\mfam$ with
extension type $\Gm$, and
We have by now obtained a complete classification of all extensions of
extensions $(M,C)$ of the family $\mfam$ such that $\xi_i \ne 0$ for
$2 \le i \le n$. Notice that all iterated extensions $(M,C)$ in this
classification is such that $M$ is indecomposable. In fact, the condition
$\xi_i \ne 0$ for $2 \le i \le n$ is equivalent with the condition that
that $M$ is indecomposable when $\mfam$ satisfy (*):

\begin{prop} \label{p:ind-nz}
Let $\mfam$ be a family satisfying (*), and let $(M,C)$ be an iterated
extension of the family $\mfam$ of length $n$ such that $M$ is
indecomposable. Then $\xi_i \ne 0$ for $2 \le i \le n$.
\end{prop}
\begin{proof}
The result is obviously true if $n \le 2$, so let us proceed by induction
on $n$: We assume that $n \ge 3$, and let $(M,C)$ be an extension of
extensions of $\mfam$ of length $n$ with $M$ indecomposable. Clearly,
$C_{n-1}$ has finite length and therefore an indecomposable decomposition
\[ C_{n-1} = N_1 \oplus \dots \oplus N_q. \]
Suppose $q > 1$. Since $N_j$ is a left $R$-module of finite length for $1
\le j \le q$, $N_j$ has a co-filtration of length $n_j < n$
\[ N_j = C_{j,n_j} \to C_{j,n_j-1} \to \dots \to C_{j,1} \to C_{j,0}=0, \]
with $K_{ji} = \ker(C_{j,i} \to C_{j,i-1}) \cong M_{l(j,i)}$ for $1 \le i
\le n_j$, with $\alpha(j,i) \in I$. Since $\xi_n \in \ext^1_R(C_{n-1},K_n)
\cong \oplus \ext^1_R(N_j,K_n)$, we can write $\xi_n = (\xi_{n,1}, \dots,
\xi_{n,q})$ with $\xi_{n,j} \in \ext^1_R(N_j,K_n)$. We claim that $\xi_{n,j}
\ne 0$ for all $j$: Assume that $\xi_{n,j}=0$ for some $j$, then we may
assume $j=1$ with no loss of generality. For each $j$, we choose a
representative $\psi(j) \in \der_k(R,\hmm_k(N_j,K_n))$ of $\xi_{n,j}$, and
let $\phi \in \hmm_k(N_1,K_n)$ satisfy $r\phi - \phi r = \psi(1)_r$ for all
$r \in R$. We may assume that $M \cong K_n \oplus (N_1 \oplus \dots \oplus
N_q)$ considered as a vector space over $k$, and that the left $R$-module
structure of $M$ is given by
\[ r(k,n_1, \dots, n_k) = (rk + \sum_{j=1}^q \psi(j)_r (n_j), rn_1, \dots,
rn_q) \]
for all $r \in R, \; (k,n_1, \dots, n_q) \in M$. Let $M' \cong M$ considered
as a vector space over $k$, and let $M'$ have a left $R$-module structure
given by
\[ r(k,n_1, \dots, n_k) = (rk + \sum_{j=2}^q \psi(j)_r (n_j), rn_1, \dots,
rn_q) \]
for all $r \in R, \; (k,n_1, \dots, n_q) \in M'$. Then the homomorphism $\pi:
M \to M'$ given by $(k,n_1, \dots, n_q) \mapsto (k+\phi(n_1),n_1, \dots,n_q)$
defines an isomorphism of left $R$-modules. But $M'$ has $N_1$ as a direct
summand, and $M'$ is indecomposable since $M$ is, so this implies $q=1$.
Since we have assumed that $q>1$, we must have $\xi_{n,j} \ne 0$ for all $j$.

Clearly, $N_j$ is indecomposable of length $n_j<n$ for $1 \le j \le q$, so by
the induction hypothesis, the extensions $\xi_{ji} \in \ext^1_R(C_{j,i-1},
K_{ji}) \ne 0$ for $2 \le i \le n_j$. So $(N_j,C_j)$ is an extension of
extensions of the family $\mfam$ which is part of the classification in
proposition \ref{p:indscl}. Let $K_n = M_{\alpha(n)}$ with $\alpha(n) \in
I$, and let $\alpha(n-1) \in I$ be the unique index such that $\alpha(n) =
\sigma(\alpha(n-1))$. From the proof of lemma \ref{l:indextnstr}, we see
that $\ext^1_R(N_j,K_n) \cong \ext^1_R(K_{j,n_j},K_n)$ and $K_{j,n_j} \cong
M_{\alpha(n-1)}$ for $1 \le j \le q$. We may therefore assume that
\begin{align*}
M &\cong K_n \oplus (N_1 \oplus \dots \oplus N_q) \\
  &\cong K_n \oplus (\oplus_{j=1}^q (\oplus_{i=1}^{n_j} K_{j,i})),
\end{align*}
considered as a vector space over $k$, and the $R$-module structure on $M$
maps $K_{j,i}$ into $K_{j,i+1}$ when $i<n_j$, and $K_{j,n_j}$ into $K_n$.

With no loss of generality, we may assume that $n_1 \le n_2$. Let us choose
representatives $\psi(j) \in \der_k(R,\hmm_k(K_{j,n_j},K_n)$ of $\tau_{n,j}
\in \ext^1_R(K_{j,n_j},K_n) \cong k^*$ for $1 \le j \le q$. Then we can
find $c \in k^*$ such that $\psi(1) = c \psi(2)$. Let us also choose
representatives $\psi(j,i) \in \der_k(R,\hmm_k(K_{j,n_j-i},K_{j,n_j-i+1}))$
corresponding to the extensions $\tau_{n_j-i+1}(N_j,C_j) \in
\ext^1_R(K_{j,n_j-i},K_{n_j-i+1})$ for $1 \le j \le 2, \; 1 \le i \le
n_j-1$. Since $K_{1,n_1} \cong K_{2,n_2}$, we also have $K_{1,n_1-i} \cong
K_{2,n_2-i}$ for $1 \le i \le n_1-1$, so we can find $c_i \in k^*$ such that
$\psi(1,i) = c_i \psi(2,i)$ for $1 \le i \le n_1-1$. Let $M' \cong M$
considered as a vector space over $k$, and let $M'$ have the left $R$-module
structure given by
\[ r(k, n_1, \dots, n_q) = (rk + \sum_{j=2}^q \psi(j)_r (k_{j,n_j}), rn_1,
\dots, r n_q) \]
for all $r \in R, \; (k, n_1, \dots, n_q) \in M'$, where we write $n_j =
(k_{j,n_j}, \dots, k_{j,1})$. Let us also write $C(i) = c_1 c_2 \dots c_i$
for $1 \le i \le n_1-1$. Then the map $\pi: M \to M'$ given by
\begin{align*}
(k, n_1, \dots, n_q) &\mapsto (k, n_1, n_2, \dots, n_q) \\
&+ (0, 0, c (k_{1,n_1}, C(1) k_{1,n_1-1}, \dots, C(n_1-1) k_{1,1}, 0, \dots,
0), 0, \dots, 0)
\end{align*}
defines an $R$-linear isomorphism of left $R$-modules. But $N_1$ is a
direct summand of $M'$, and $M'$ is indecomposable since $M$ is, so this
implies $q=1$. We must therefore conclude that $C_{n-1}$ is indecomposable.
But this implies that $\xi_{n-1} \ne 0$. By induction, it follows that
$\xi_i \ne 0$ for $2 \le i \le n$.
\end{proof}

\begin{thm} \label{t:exclass}
Let $\mfam$ be a family satisfying (*), and let $\Gm$ be an ordered quiver.
There exists an iterated extension of the family $\mfam$ with extension
type $\Gm$ such that $M$ is indecomposable if and only if $\Gm = \Gm(\alpha,n)$
for some admissible couple $(\alpha,n)$ in $I \times \n$. Moreover,
$\indextncl(\mfam,\Gm) \cong (k^*)^{n-1}$ and $\indmodfcl(\mfam,\Gm) = \{ * \}$
in this case.
\end{thm}
\begin{proof}
Proposition \ref{p:ind-nz} shows that $\indextncl(\mfam,\Gm) = \nsextncl(\mfam,
\Gm)$ when $\mfam$ satisfy (*). The rest is clear.
\end{proof}

Let us denote by $M(\alpha,n)$ the $R$-module representing the isomorphism
class $*$ in the above classification for each admissible couple $(\alpha,n)
\in I \times \n$. We remark that the proofs given in this section are
constructible, and therefore we may in principle construct $M(\alpha,n)$. In
fact, the simple modules $\mfam$ and their extensions is enough to construct
the left $R$-module $M(\alpha,n)$.

\begin{thm}
Let $k$ be an algebraically closed field, let $R$ be an associative
$k$-algebra, and let $\mfam = \{ M_{\alpha}: \alpha \in I \}$ be a family
of orthogonal points in $\rmod$ such that the species of $\mfam$ is a
$k$-quiver. Then the category $\modf$ is uniserial if and only if $\mfam$
satisfies the condition (*). If this is the case, there is a complete
classification of all indecomposable modules in $\modf$.
\end{thm}
\begin{proof}
The modules $M(\alpha,n)$ are uniserial when $(\alpha,n) \in I \times \n$
is an admissible couple and $\mfam$ satisfy (*). So if $\mfam$ satisfy (*),
all indecomposable modules in $\modf$ are uniserial. Conversely, if all
indecomposable modules in $\modf$ are uniserial, then $\mfam$ satisfy (*)
by lemma \ref{l:ce1}, \ref{l:ce2} and \ref{l:ce3}.
\end{proof}

In the uniserial case, the list of indecomposable modules in $\modf$ is
of course given by
\[ \{ M(\alpha,n): (\alpha,n) \in I \times \n \text{ is an admissible
couple } \}. \] Moreover, we have seen that the left $R$-modules
$M(\alpha,n)$ are constructible. That is, given the simple modules
in $\mfam$ and the extensions $\tau_i \in \ext^1_R(K_{i-1},K_i)$
expressed in Hochschild cohomology, we can construct the modules
$M(\alpha,n)$ in concrete terms (for instance in terms of
generators and relations).

\section{Applications: Regular holonomic D-modules in dimension 1}

We find many examples of uniserial length categories among the 
categories of regular holonomic D-modules on curves. We shall 
prove this in a number of examples, and at the same time
give the corresponding classification of all indecomposable 
objects.  

In most cases, this gives new proofs of known results. But the 
classification of graded holonomic $D$-modules when $D$ is the 
first Weyl algebra $D=\weyl$ or the ring of differential 
operators $D = \diff(A)$ over an affine monomial curve $A$ is 
new. 

\subsection{The local analytic case}

Let $k = \cc$ be the complex numbers, and consider the local ring
$A = k \{ t \}$ of convergent power series with coefficients in
$k$. Let furthermore $D = \diff(A) = A <\partial>$ be the ring of
$k$-linear differential operators on $A$, with $\partial = d/dt$.
Explicitly, $D$ is the skew polynomial ring
\begin{equation*}
D = \{ \; \sum_{i=0}^d p_i \partial^i : \; d \ge 0, \; p_i \in A
\text{ for } 0 \le i \le d \; \}
\end{equation*}
with the relation $\partial t = t \partial + 1$.

Let us consider the holonomic $D$-modules with regular
singularities. These were completely classified by Boutet de
Monvel \cite{bou83}. The same classification result was later
obtained by Brian\c{c}on, Maisonobe \cite{br-ma84, br-ma86}, using
division algorithms and perverse sheaves. For definitions of the
terms \emph{holonomic} $D$-modules with \emph{regular
singularities}, we refer the reader to van den Essen \cite{es87}.
However, for the ring $D$ considered in this section, it is useful
to notice the following facts: A $D$-module is holonomic if and
only if it has finite length. Moreover, if $0 \to M' \to M \to M''
\to 0$ is a short exact sequence of holonomic $D$-modules, then
$M$ has regular singularities if and only if $M'$ and $M''$ has
regular singularities.

Let $I = \{ \alpha \in k: 0 < \re(\alpha) < 1 \} \cup \{ 0,1 \}$,
and let $\mfam = \{ M_{\alpha}: \alpha \in I \}$ be the family of
left $R$-modules given by $M_0 = D/D\partial, \; M_1 = D/Dt$, and
$M_{\alpha} = D/D(E-\alpha)$ with $E = t\partial$ for all $\alpha
\in I \setminus \{ 0,1 \}$. It is well-known that $\mfam$ is the
family of all simple left $D$-modules with regular singularities,
up to isomorphism. From the comments above, we see that the
category $\modf$ is the category of regular holonomic $D$-modules.

Since $\mfam$ consists of simple, non-isomorphic $D$-modules, it
is clear that $\mfam$ is a family of orthogonal points in
$\modc(D)$. Moreover, it is easy to check that the species of
$\mfam$ is a $k$-quiver satisfying the condition (*) of section
\ref{s:uniserial}. It follows that the category $\modf$ of regular
holonomic $D$-modules is a uniserial category. Moreover, each
indecomposable object $M$ in $\modf$ can be constructed explicitly
with the methods from section \ref{s:uniserial}. With our methods,
we therefore reprove the classification of regular holonomic
$D$-modules:

\begin{thm} \label{t:localana}
Let $k = \cc$ and let $D$ be the ring of differential operators on
the $k$-algebra $A = k \{ t \}$ of convergent power series with
coefficients in $k$. Then the category of regular holonomic 
$D$-modules is uniserial. Moreover, any regular holonomic $D$-module 
is a finite direct sum of the indecomposable ones, given by
\[ \{ M(\alpha,n): \alpha \in I, \; n \ge 1 \}, \]
where $M(\alpha,n) = D / D \; w(\alpha,n)$ for $\alpha = 0,1$, and
$M(\alpha,n) = D / D (E - \alpha)^n$ for $\alpha \in I \setminus 
\{ 0,1 \}$.
\end{thm}

For $\alpha = 0,1$, we use the notation $w(\alpha,n)$ for the
alternating word in the letters $t, \partial$ of length $n$,
ending in $\partial$ if $\alpha = 0$ and ending in $t$ if $\alpha
= 1$. Notice that any couple $(\alpha,n) \in I \times \n$ is
admissible.

Let $A' \subseteq A$ be a sub-algebra of $A = k \{ t \}$ such that
$\dim_k A/A'$ is finite. Then $D' = \diff(A')$ is Morita
equivalent to $D = \diff(A)$ by Smith, Stafford \cite{sm-st88},
proposition 3.3, and this Morita equivalence preserves regular
holonomic modules by van den Essen \cite{es89}, theorem 3.1. It
follows that the category of regular holonomic $D'$-modules is
uniserial, and there is a classification of regular holonomic
$D'$-modules similar to theorem \ref{t:localana}.

Let $A_n = k \{ t_1, \dots, t_n \}$ be the ring of convergent
power series in $n$ variables with coefficients in $k$, and let
$D_n = \diff(A_n)$ be the ring of differential operators on $A_n$.
The category of regular holonomic $D_n$-modules supported by an
irreducible analytic curve is equivalent to the category of
regular holonomic $D$-modules (with $D = D_1$ as above), by van
Doorn, van den Essen \cite{do-es87}. It follows that the category
of regular holonomic $D_n$-modules supported by an irreducible
analytic curve is uniserial, and there is a classification of
regular holonomic $D$-modules supported by an irreducible analytic
curve similar to theorem \ref{t:localana}.

\subsection{The formal case}

Let $k$ be an algebraically closed field of characteristic $0$,
and let $A = k[[t]]$ be the local ring of formal power series
in $t$. Let furthermore $D = \diff(A) = A < \partial > = \weylf$ 
be the ring of $k$-linear differential operators on $A$, with 
$\partial = d/dt$. Explicitly, $D$ is the skew polynomial ring
\[ D = \{ \; \sum_{i=0}^d p_i \partial^i : \; d \ge 0, \; p_i
\in A \text{ for } 0 \le i \le d \; \} \]
with the relation $\partial t = t \partial + 1$.

Let us consider the holonomic $D$-modules. These fall into $2$ 
classes: The regular holonomic $D$-modules and the irregular 
holonomic $D$-modules. The classification of the regular ones 
is completely parallell to the classification of regular 
holonomic $D$-modules in the local analytic case. Again, we 
refer the reader to van den Essen \cite{es87} for definitions 
of the terms holonomic and regular holonomic $D$-modules. For 
the ring $D$ we consider in this section, it remains valid that 
a module is holonomic if and only if it has finite length, and 
for a short exact sequence of holonomic $D$-modules $0 \to M' 
\to M \to M'' \to 0$, we have that $M$ is regular if and only 
if $M'$ and $M''$ are regular.

The classification of the irregular holonomic $D$-modules were
first obtained by Puninski \cite{pu00}. It is based upon the 
classification of the simple irregular $D$-modules of van den
Essen, Levelt \cite{es-le92}. We remark that this classification
is very similar to the classification of simple modules over 
the first Weyl algebra in Block \cite{bl81}: The simple 
$D$-modules which are torsion modules over the sub-module $A 
\subseteq D$ is given by $\{ D/Dt \}$, up to isomorphism. The
simple $D$-modules which are torsion free $A$-modules are 
parametrized by the simple $\diff(K)$-modules, where $K=k((t))$
is the quotient field of $A$. 

It is clear that any simple $\diff(K)$-module is of the form
$N = \diff(K) / \diff(K) P$, where $P \in \diff(K)$ is 
irreducible, since $\diff(K)$ is a left and right principal
ideal domain. Moreover, the simple $A$-torsion free module $M$ 
corresponding to $N$ is given by $\soc_D N$. But if $N \cong K$,
then $\soc_D N = A$, and otherwise $\soc_D N = N$. We denote by
$I'$ the set of equivalence classes of irreducible elements in
$\diff(K)$, and by $M_{\alpha}$ the simple, $A$-torsion free 
$D$-module corresponding to $\alpha \in I'$. 

Let $I = I' \cup \{ t \}$, and let $M_t = D/Dt$ be the simple 
$D$-module which has $A$-torsion. Then $\mfam = \{ M_{\alpha}:
\alpha \in I \}$ is the set of simple $D$-modules (up to 
isomorphism) by van den Essen, Levelt \cite{es-le92}. This is 
clearly a family of orthogonal points in $\modc(D)$, since it 
consists of simple, pairwise non-isomorphic $D$-modules. 
Moreover, it follows from Puninski \cite{pu00}, fact 2.3 that 
the species of $\mfam$ is a $k$-quiver, and from Puninski 
\cite{pu00}, proposition 3.4 that this species satisfy condition 
(*) of section \ref{s:uniserial}. 

If follows that the category of holonomic $D$-modules is 
uniserial, and we can construct the holonomic $D$-modules 
explicitly with the methods of section \ref{s:uniserial}.
This gives a new proof of Puninskis classification of holonomic
$D$-modules:

\begin{thm}[Puninski] \label{t:formal}
Let $k$ be a algebraically closed field of characteristic $0$,
and let $D$ be the ring of differential operators on the ring 
$A=k[[t]]$ of formal power series with coefficients in $k$. 
Then the category of holonomic $D$-modules is uniserial. 
Moreover, any holonomic $D$-module is a finite direct sum of 
the indecomposable ones, given by 
\[ \{ M(\alpha,n): \alpha \in I, \; n \ge 0 \}, \]
where $M(\alpha,n) = D/D w(\alpha,n)$ for $\alpha = t, \partial$
and $M(\alpha,n) = M_{n-1}(\alpha)$ in Puninskis notation for 
$\alpha \in I \setminus \{ t, \partial \}$. 
\end{thm}

For $\alpha = t, \partial$, we use the notations $w(\alpha,n)$ 
for the alternating word in the letter $t, \partial$ of length
$n$ ending in $\alpha$. Notice that all couples $(\alpha,n) \in 
I \times \n$ are admissible. We also remark that it does not 
seem to be possible to use the same methods to obtain a complete 
classification of holonomic $D$-module in the local analytic 
case.

Let $D_n = \diff(k[[t_1, \dots, t_n]])$, and consider the 
category of holonomic $D_n$-modules supported by an irreducible 
curve. By the equivalence of categories between holonomic 
$D_n$-modules supported on an irreducible curve and holonomic 
$D_1$-modules, we conclude that the category of holonomic 
$D_n$-modules supported by an irreducible curve is uniserial, 
and a classification result similar to theorem \ref{t:formal} 
holds for this category.

\subsection{The graded algebraic case}

Let $k$ be an algebraically closed field of characteristic $0$,
and let $A = k[t]$ be the polynomial ring in one variable over
$k$. Let furthermore $D = \diff(A) = A < \partial > = \weyl$ be
the first Weyl algebra, with $\partial = d/dt$. Explicitly, $D$ 
is the skew polynomial ring
    \[ D = \{ \; \sum_{i=0}^d p_i \partial^i : \; d \ge 0, \; p_i 
    \in A \text{ for } 0 \le i \le d \; \}, \]
with the relation $\partial t = t \partial + 1$.

The ring $A = k[t]$ is a $\z$-graded $k$-algebra in a natural way,
such that $t^i$ is homogeneous of degree $i$ for all integers $i
\ge 0$. This induces a natural $\z$-grading of the Weyl algebra
$D$, such that $D$ is a $\z$-graded $k$-algebra: We say that a
differential operator $P \in D$ is homogeneous of weight $w \in
\z$ if $P*A_i \subseteq A_{i+w}$ for all integers $i \in \z$.
Explicitly, $D$ is the $\z$-graded $k$-algebra generated by the
homogeneous monomials $t^i \partial^j$ of weight $i-j$ for all
$i,j \ge 0$.

Let $\ccat = \grdmod$ be the category of graded $D$-modules. The 
objects of $\grdmod$ is the $\z$-graded, left $D$-modules, and the 
morphisms are the homogeneous homomorphisms (of any degree) between 
graded $D$-modules. In this section, we shall study the full
sub-category $\grhol(D)$ of $\grdmod$, consisting of all objects in 
$\grdmod$ which have finite length. Even though the notations we 
have introduced and the results we have obtained in this paper are 
stated for length categories in $\ccat = \rmod$, we shall feel free
to use them in the case $\ccat = \grdmod$ as well, see the note in 
the introduction.

Let $M$ be a graded $D$-module. We shall denote by $\vect M$ the 
module $M$ considered as a $D$-module, forgetting the graded 
structure. First, notice that $M$ is simple in the category 
$\grdmod$ if and only if $\vect M$ is simple in $\modc(D)$: One 
implication is obvious, the other follows from N\u{a}st\u{a}sescu, 
van Oystaeyen \cite{na-oy79}, theorem II.7.5 and the fact that if 
$M$ is a simple object of $\grdmod$, then $\vect M$ is a $D$-module
of finite length. 

It follows that a graded $D$-module $M$ is of finite length if and
only if $\vect M$ is a $D$-module of finite length. On the other 
hand, it is well known that for the first Weyl algebra $D=\weyl$, a 
$D$-module is holonomic if and only if it has finite length. So a 
graded $D$-module has finite length if and only if $\vect M$ is 
holonomic. It is therefore natural to define the \emph{category of 
graded holonomic $D$-modules} to be the category $\grhol(D)$. We 
notice that if $M$ is an object of $\grhol(D)$, then $\vect M$ is a 
regular holonomic $D$-module in the sense of van den Essen 
\cite{es87}.

Let us classify the simple objects in $\grhol(D)$: We know that these
are exactly the simle graded $D$-modules. On the other hand, all 
simple $D$-modules have been classified by Block \cite{bl81}. We 
shall adapt this classification to obtain a classification of all
simple objects in $\grhol(D)$, up to graded isomorphism.

Let $M$ be a simple graded $D$-module. Then $\vect M$ is either 
torsion or torsion free considered as a module over the sub-ring 
$A \subseteq D$. If $\vect M$ is a torsion module over $A$, then
$\vect M \cong D/D(t-\alpha)$ for some $\alpha \in k$ by Block
\cite{bl81}, proposition 4.1 and corollary 4.1. So $M \cong D/Dt$
in $\grhol(D)$, where $D/Dt$ has the natural graded structure 
inherited from $D$. 

There is a bijective correspondence, given by localization, between
simple graded $D$-modules which are torsion free over $A$ and simple 
graded $\diff(T)$-modules, where $T=k[t,t^{-1}]$. This follows from 
Block \cite{bl81}, lemma 2.2.1 and corollary 2.2, slightly adapted
to the graded situation. But any simple graded $\diff(T)$-module is
of the form $N_{\alpha} = \diff(T) / \diff(T) (E - \alpha)$ for some
$\alpha \in k$, where $E = t \partial$. Moreover, $N_{\alpha} \cong
N_{\beta}$ as graded modules if and only if $\alpha - \beta \in \z$.

Let $I' \subseteq k$ be a subset of $k$ containing $0$ such that 
the natural composition $I' \to k \to k / \z$ is a bijection. For 
each non-zero $\alpha \in I'$, the simple graded $D$-module 
$M_{\alpha} = D/D(E-\alpha)$ corresponds to $N_{\alpha}$ in the 
correspondence above. Moreover, the simple graded $D$-module $M_0 = 
D/D \partial$ corresponds to $N_0$. 

Let $I = I' \cup \{ 1 \}$, and let $M_1 = D/Dt$. It follows that 
$\mfam = \{ M_{\alpha}: \alpha \in I \}$ is a family of simple graded
$D$-modules with the following property: Any simple graded $D$-module 
$M$ is isomorphic to $M_{\alpha}$ (as a graded $D$-module) for a 
unique $\alpha \in I$. In this sense, $\mfam$ is the family of simple
objects in $\grhol(D)$, up to graded isomorphism.

Clearly, $\mfam$ is a family of orthogonal points in $\grdmod$, since
it is a family of simple, non-isomorphic objects. It is also easy to
check that the species of $\mfam$ is a $k$-quiver which satisfy the 
condition (*) of section \ref{s:uniserial}. So the sub-category 
$\modf$ of $\grdmod$ is a uniserial category. In fact, $\modf$ is the
length category $\grhol(D)$ of graded holonomic $D$-modules. So we 
conclude that the category $\grhol(D)$ of graded holonomic $D$-modules 
is a uniserial category. Moreover, each indecomposable object in 
$\grhol(D)$ can be constructed explicitly with the methods of section
\ref{s:uniserial}. We obtain the following classification result: 

\begin{thm} \label{t:grhol}
Let $k$ be an algebraically closed field of characteristic $0$, and 
let $D = \weyl$ be the first Weyl algebra over $k$. Then the category 
of graded holonomic $D$-modules is uniserial. Moreover, any graded 
holonomic $D$-module is a finite direct sum of the indecomposable 
ones, given by
\[ \{ M(\alpha,n): \alpha \in I, \; n \ge 1 \}, \]
where $M(\alpha,n) = D/D w(\alpha,n)$ for $\alpha = 0,1$ and $M(\alpha,
n) = D/D(E-\alpha)^n$ for $\alpha \in I \setminus \{ 0,1 \}$. 
\end{thm}

For $\alpha = 0,1$, we use the notation $w(\alpha,n)$ for the 
alternating word in the letters $t,\partial$ of length $n$, ending in 
$\partial$ if $\alpha = 0$ and ending in $t$ in $\alpha = 1$. Notice 
that any couple $(\alpha,n) \in I \times \n$ is admissable.

Let $A'$ be an affine monomial curve over $k$, and let $D' = \diff(A')$
be the corresponding ring of differential operators. By Smith, Stafford
\cite{sm-st88}, the category of holonomic $D$-modules and the category
of holonomic $D'$-modules are equivalent, and clearly graded structures
are conserved under this equivalence. It follows that the category of 
graded holonomic $D'$-modules is uniserial, and there is a classification
of graded holonomic $D'$-modules similar to theorem \ref{t:grhol} for any
affine monomial curve $A'$.

\section{The wild case}

Let $k$ be an algebraically closed field, let $R$ be an associative
$k$-algebra, and let $\mfam$ be a family of orthogonal points in 
$\rmod$ such that $\enm_R(M_{\alpha}) = k$ for all $\alpha \in I$.
In this section, we shall mention a result which gives a sufficient
condition for the length category $\modf$ to be wild in a strong 
sense.

Let $W = k <x,y>$ be the free associative $k$-algebra on two 
generators, and consider the category $\modw$ of left $W$-modules 
which are finitely dimensional as vector spaces over $k$. We say 
that the category $\modf$ is \emph{wild} if there is a full exact
embedding of $\modw$ into $\modf$, following Klingler, Levy 
\cite{kl-le95}. It is well-known that if $\modf$ is wild in this
sense, a classification of the indecomposable modules in $\modf$
would contain a classification of all indecomposable modules of 
finite dimension over $k$ over any finite dimensional $k$-algebra.

In Klingler, Levy \cite{kl-le95}, it was shown that the category 
of holonomic modules over the first Weyl algebra is wild in the 
above sense if $k$ has characteristic $0$. In fact, a full exact 
embedding can be chosen such that every module in its image has 
socle-height $2$.

\begin{thm}[Klingler-Levy]
Let $k$ be an algebraically closed field, let $R$ be an associative
$k$-algebra, and let $\mfam$ be a family of orthogonal points in
$\rmod$ such that $\enm_R(M_{\alpha}) = k$ for all $\alpha \in I$.
If the Gabriel quiver of $\modf$ contains the quiver $Q_5$ given
by 
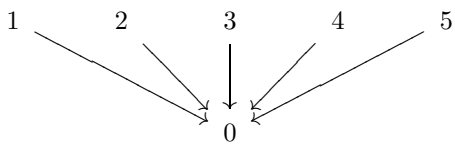
\[ \xymatrix{ 1 \ar[drr] & 2 \ar[dr] & 3 \ar[d] & 4 \ar[dl] & 5 \ar[dll] \\
& & 0 & & } \] 
then there is a full exact embedding of $\modw$ into $\modf$, and 
all modules in the image of this embedding has socle-height $2$. In 
particular, $\modf$ is wild in this case.
\end{thm}
\begin{proof}
This follows from Klingler, Levy \cite{kl-le95}, theorem 2.12 with 
some minor changes.
\end{proof}

\bibliographystyle{amsplain}
\bibliography{main}

\providecommand{\bysame}{\leavevmode\hbox to3em{\hrulefill}\thinspace}
\begin{thebibliography}{10}

\bibitem{am-ri68}
Ivar~Kr. Amdal and Frode Ringdal, \emph{Cat\'egories unis\'erielles}, Comptes
  Rendus de l'Academie des Sciences \textbf{267} (1968), A85--87, A247--249.

\bibitem{bl81}
Richard~E. Block, \emph{The irreducible representations of the lie algebra
  $sl_2$ and of the {W}eyl algebra}, Advances in Mathematics \textbf{39}
  (1981), 69--110.

\bibitem{bou83}
Louis Boutet~de Monvel, \emph{D-modules holon\^omes r\'eguliers en une
  variable}, Math\'ematique et Physique, S\'eminaire de l'E.N.S. 1979-82,
  Progress in Mathematics, vol.~37, Birkh\"auser, 1983, pp.~313--321.

\bibitem{br-ma84}
J.~Brian\c{c}on and Ph. Maisonobe, \emph{Id\'eaux de germes d'op\'erateurs
  diff\'erentiels \`a une variable}, L'Enseignement Math\'ematique \textbf{30}
  (1984), 7--38.

\bibitem{br-ma86}
\bysame, \emph{Germes de {D}-modules \`a une variable et leurs solutions},
  Preprint~88, Universit\'e de Nice, January 1986.

\bibitem{de-xi98}
Bangming Deng and Jie Xiao, \emph{A quiver description of hereditary categories
  and its application to the first {W}eyl algebra}, Algebras and modules {II}
  (Idun Reiten, Sverre~O. Smal\"o, and \"Oyvind Solberg, eds.), Candian
  Mathematical Society Conference Proceedings, vol.~24, American Mathematical
  Society, 1998, pp.~125--137.

\bibitem{er02}
Eivind Eriksen, \emph{Non-commutative deformations of modules (after
  {L}audal)}, Preprint no. 1, University of Warwick, 2002.

\bibitem{ga73}
Peter Gabriel, \emph{Indecomposable representations {II}}, Symposia
  Mathematica, vol.~XI, Academic Press, 1973, pp.~81--104.

\bibitem{ka83}
Victor~G. Kac, \emph{Root systems, representations of quivers and invariant
  theory}, Invariant theory, Springer Lecture Notes, vol. 996, Springer, 1983,
  pp.~74--108.

\bibitem{kl-le95}
Lee Klingler and Lawrence~S. Levy, \emph{Wild torsion modules over {W}eyl
  algebras and general torsion modules over {HNP}s}, Journal of Algebra
  \textbf{172} (1995), 273--300.

\bibitem{lau95}
Olav~Arnfinn Laudal, \emph{Non-commutative deformations of modules}, Preprint
  no. 2, University of Oslo, 1995.

\bibitem{lau98}
\bysame, \emph{Non-commutative deformations of modules {II}}, Unpublished
  manuscript, 1998.

\bibitem{lau00}
\bysame, \emph{Noncommutative algebraic geometry}, Preprint no. 115, Max Planck
  Institute of Mathematics, 2000.

\bibitem{mc-ro73}
J.C. McConnell and J.C. Robson, \emph{Homomorphisms and extensions of modules
  over certain differential polynomial rings}, Journal of Algebra \textbf{26}
  (1973), 319--342.

\bibitem{na-oy79}
C.~N\u{a}st\u{a}sescu and F.~Van~Oystaeyen, \emph{Graded and filtered rings and
  modules}, Lecture notes in mathematics, no. 758, Springer--Verlag, 1979.

\bibitem{pu00}
Gennadi Puninski, \emph{Finite lenght and pure-injective modules over a ring of
  differential operators}, Journal of Algebra \textbf{231} (2000), 546--560.

\bibitem{qu69}
Daniel Quillen, \emph{On the endormorphism ring of a simple module over an
  enveloping algebra}, Proceedings of the American Mathematical Society
  \textbf{21} (1969), 171--172.

\bibitem{ri76}
Claus~Michael Ringel, \emph{Representations of {K}-species and bimodules},
  Journal of Algebra \textbf{41} (1976), 269--302.

\bibitem{sm-st88}
S.~Paul Smith and J.~T. Stafford, \emph{Differential operators on an affine
  curve}, Proceedings of the London Mathematical Society. Third Series
  \textbf{56} (1988), no.~2, 229--259.

\bibitem{es87}
Arno van~den Essen, \emph{Modules with regular singularities over filtered
  rings and algebraic micro-localization}, S\'eminaire d'Alg\`ebre Paul Dubreil
  et Marie-Paule Malliavin, Lecture notes in mathemtics, vol. 1296, Springer,
  1987, pp.~125--1570.

\bibitem{es89}
\bysame, \emph{Modules with regular singularities on a curve}, Journal of the
  London Mathematical Society. Second Series \textbf{40} (1989), 193--205.

\bibitem{es-le92}
Arno van~den Essen and A.~Levelt, \emph{An explicit description of all simple
  $k[[x]][\partial]$-modules}, Contemporary Mathematics, vol. 130, American
  Mathematical Society, 1992, pp.~121--131.

\bibitem{do-es87}
M.G.M. van Doorn and Arno van~den Essen, \emph{D-modules with support on a
  curve}, Publ. Res. Inst. Math. Sci. \textbf{23} (1987), 937--953.

\end{thebibliography}

\end{document}